\documentclass[a4paper, 10pt]{article}

\usepackage{verbatim}
\usepackage{amsmath}
\usepackage{amscd}
\usepackage{amssymb}
\usepackage{amsthm}
\usepackage{graphicx}
\usepackage[all]{xy}

\newtheorem{thm}{Theorem}[section]
\newtheorem{theorem}[thm]{Theorem}
\newtheorem{prop}[thm]{Proposition}
\newtheorem{proposition}[thm]{Proposition}
\newtheorem{cor}[thm]{Corollary}
\newtheorem{corollary}[thm]{Corollary}
\newtheorem{lem}[thm]{Lemma}
\newtheorem{lemma}[thm]{Lemma}
\newtheorem{step}{Step}

\newtheorem{claim}{Claim}

\theoremstyle{definition}
\newtheorem{df}[thm]{Definition}
\newtheorem{definition}[thm]{Definition}
\newtheorem{rem}[thm]{Remark}

\newtheorem{expl}[thm]{Example}

\newcommand{\dge}{\rotatebox[origin=c]{45}{$\ge$}}
\newcommand{\uge}{\rotatebox[origin=c]{315}{$\ge$}}
\newcommand{\deq}{\rotatebox[origin=c]{45}{$=$}}
\newcommand{\ueq}{\rotatebox[origin=c]{315}{$=$}}
\newcommand{\updots}{\hbox to1.65em{\rotatebox[origin=c]{45}{$\cdots$}}}
\newcommand{\dndots}{\hbox to1.65em{\rotatebox[origin=c]{315}{$\cdots$}}}
\newcommand{\rueq}{\rotatebox[origin=c]{45}{$=$}}
\newcommand{\lueq}{\rotatebox[origin=c]{315}{$=$}}
\newcommand{\lmd}[1]{\hbox to1.65em{$\hfill \lambda_{#1} \hfill$}}
\newcommand{\llmd}[2]{\hbox to1.65em{$ \lambda_{#2}^{(#1)}$}}

\newcommand{\tr}{\mathrm{tr}\,}
\newcommand{\re}{\Re}
\newcommand{\im}{\Im}

\newcommand{\simto}{\xrightarrow{\sim}}

\newcommand{\fanr}{r}

\newcommand{\ltilde}{\widetilde{l}}

\newcommand{\Xtilde}{\widetilde{X}}

\newcommand{\varphitilde}{\widetilde{\varphi}}
\newcommand{\Phitilde}{\widetilde{\Phi}}

\newcommand{\omegatilde}{\widetilde{\omega}}
\newcommand{\virt}{\mathrm{virt}}
\newcommand{\scMbar}{\overline{\scM}}
\newcommand{\Cr}{\operatorname{Cr}}
\newcommand{\scD}{\mathcal{D}}
\newcommand{\scV}{\mathcal{V}}
\newcommand{\nofaces}{m}

\newcommand{\be}{\mathbf{e}}
\newcommand{\boldf}{\mathbf{f}}
\newcommand{\frakL}{\mathfrak{L}}
\newcommand{\frakP}{\mathfrak{P}}
\newcommand{\fraku}{\mathfrak{u}}
\newcommand{\frakv}{\mathfrak{v}}
\newcommand{\frakx}{\mathfrak{x}}
\newcommand{\fraky}{\mathfrak{y}}

\newcommand{\Ker}{\operatorname{Ker}}
\newcommand{\Image}{\operatorname{Im}}

\newcommand{\bC}{\mathbb{C}}
\newcommand{\bF}{\mathbb{F}}
\newcommand{\bN}{\mathbb{N}}
\newcommand{\bP}{\mathbb{P}}
\newcommand{\bQ}{\mathbb{Q}}
\newcommand{\bR}{\mathbb{R}}
\newcommand{\bZ}{\mathbb{Z}}
\newcommand{\frakm}{\mathfrak{m}}
\newcommand{\frakX}{\mathfrak{X}}

\newcommand{\scL}{\mathcal{L}}
\newcommand{\scM}{\mathcal{M}}
\newcommand{\scO}{\mathcal{O}}

\newcommand{\bCx}{\bC^\times}

\newcommand{\Ad}{\operatorname{Ad}}
\newcommand{\diag}{\operatorname{diag}}
\newcommand{\Vol}{\operatorname{Vol}}
\newcommand{\Spec}{\operatorname{Spec}}
\newcommand{\PD}{\operatorname{PD}}
\newcommand{\Int}{\operatorname{Int}}
\newcommand{\ev}{\mathop{ev}\nolimits}
\newcommand{\vdim}{\operatorname{virt.dim}}

\newcommand{\Gr}{\mathop{Gr}}

\newcommand{\po}{\mathfrak{PO}}

\newcommand{\scMwhat}{\widehat{\scM}_{\mathrm{weak}}}

\newcommand{\dimh}{h} 

\title{Toric degenerations of Gelfand-Cetlin systems
and potential functions}
\author{Takeo Nishinou, Yuichi Nohara, and Kazushi Ueda}
\date{}

\begin{document}

\maketitle

\begin{abstract}
We define a toric degeneration of an integrable system
on a projective manifold,
and prove the existence of a toric degeneration
of the Gelfand-Cetlin system on the flag manifold
of type $A$.
As an application,
we calculate the potential function
for a Lagrangian torus fiber
of the Gelfand-Cetlin system.
\end{abstract}

\section{Introduction}
It is well known that a polarized toric variety $(X,\mathcal{L})$
is related to a convex polytope $\Delta_{\scL}$,
the moment polytope,
in two different ways:
\begin{itemize}
  \item $\Delta_{\scL}$ is the image of the moment map for the standard 
        torus action on $X$, and
  \item the space $H^0(X,\scL)$ of holomorphic sections of $\scL$ 
        has a basis consisting of Laurent monomials, or equivalently, 
        the weight decomposition of $H^0(X,\scL)$ with respect to 
        the torus action is multiplicity-free, and each monomial corresponds to 
        an integral point of $\Delta_{\scL}$.
\end{itemize}
Similar relations are known also for flag manifolds.
Let
$$
 \lambda = (\lambda_1 \ge \lambda_2 \ge \dots \ge \lambda_{n})
$$
be a non-increasing sequence of real numbers
and consider the orbit
$$
 \scO_\lambda = \Ad_{U(n)} \cdot \diag(\lambda_1, \dots, \lambda_n)
$$
of the Hermitian matrix
$\diag(\lambda_1, \dots, \lambda_n)$
under the adjoint action of the unitary group $U(n)$.
This orbit has a natural K\"{a}hler structure,
where the complex structure comes
from an identification with
the flag manifold $F=GL(n, \bC)/P$ of type $A$,
and the K\"{a}hler form $\omega_{\lambda}$ comes
from the Kostant-Kirillov symplectic form
on the coadjoint orbit in the dual space
$\fraku(n)^*$
of the Lie algebra $\fraku(n)$,
identified with $\scO_{\lambda}$
by the Killing form.

When all $\lambda_i$ are integral,
there is a $U(n)$-equivariant ample line bundle $\mathcal{L}_{\lambda}$
on $\scO_{\lambda}$
whose first Chern class $c_1(\mathcal{L}_{\lambda})$
is represented by $\omega_{\lambda}$.
The Borel-Weil theory states that
$H^0(F, \mathcal{L}_{\lambda})$ is an irreducible representation
of $U(n)$ of highest weight $\lambda$.
In this setting, a convex polytope $\Delta_{\lambda}$, 
called the {\it Gelfand-Cetlin polytope},
appears in two different ways:
\begin{itemize}
 \item
$\Delta_\lambda$ is the image of the moment map
of a completely integrable system on $F$
called the {\em Gelfand-Cetlin system}
\cite{GS}.
 \item
$H^0(F, \scL_{\lambda})$
admits a multiplicity-free decomposition
into one-dimensional subspaces
with respect to the action of a chain
$$
 U(1) \subset U(2) \subset \dots \subset U(n-1) \subset U(n)
$$
of subgroups.
Each of these subspaces is parametrized
by a sequence
$$
 (\lambda^{(1)}, \lambda^{(2)}, \dots, \lambda^{(n-1)}, \lambda^{(n)})
$$
of highest weights,
which is in one-to-one correspondence
with an integral point of $\Delta_{\lambda}$.
By choosing a non-zero element of each subspace,
one obtains the {\em Gelfand-Cetlin basis}
of $H^0(F, \scL_{\lambda})$
\cite{GC}.
\end{itemize}

Despite the similarities
between the toric moment map
and the Gelfand-Cetlin system,
there are marked differences:
The torus action on a flag manifold
induced by the Gelfand-Cetlin system
does not preserve the complex structure
unlike the case of a toric variety.
Although the fibers
over the interior of the moment polytope
are Lagrangian tori
in both cases,
the fibers over the boundary
of $\Delta_\lambda$ are not necessarily isotropic tori
in contrast to the case of $\Delta_{\scL}$.
In Example \ref{eg:F3},
we will see that
the moment polytope
for the full flag manifold in dimension three
has a vertex
where the fiber is a Lagrangian three-sphere.

It is known \cite{GL, B, KM}
that there is a flat family
$
 f : (\frakX, \frakL) \to \bC
$
of polarized varieties
such that
\begin{itemize}
 \item
$
 (X_t = f^{-1}(t), \scL_t = \frakL|_{X_t})
$
is isomorphic to $(F, \scL_{\lambda})$
as a polarized manifold for any $t \ne 0$,
and
 \item
$
 (X_0, \scL_0)
$
is the polarized toric variety
associated with the Gelfand-Cetlin polytope $\Delta_{\lambda}$.
\end{itemize}
In this paper,
we study the relation
between the Gelfand-Cetlin system on $(F, \scL_{\lambda})$
and the moment map of $(X_0, \scL_0)$.
To state our main result,
we make the following definition:
\begin{df} \label{def:toric_degeneration}
Let $(X, \omega)$ be a projective manifold $X$
with a K\"{a}hler form $\omega$ and
$
 \Phi : X \to \bR^N
$
be a completely integrable system
on it.
A {\em toric degeneration} of $\Phi$
consists of a flat family
$
 f : \frakX \to B
$
of algebraic varieties over a complex manifold $B$,
a K\"{a}hler form $\omegatilde$ on $\frakX$,
a piecewise smooth path $\gamma : [0,1] \to B$,
a continuous map
$
 \Phitilde : \frakX |_{\gamma([0,1])} \to \bR^N
$
on the total space 
$\frakX |_{\gamma([0,1])} = f^{-1}(\gamma([0,1]))$ 
of the family restricted to the path,
and a flow $\phi_t$ on $\frakX |_{\gamma([0,1])}$ which 
covers the path $\gamma$
such that
\begin{itemize}
 \item
$\Phi_t = \Phitilde|_{X_t}$ is a completely integrable system
on the K\"{a}hler variety \\
$(X_t = f^{-1}(\gamma(t)), \omega_t = \omegatilde|_{X_t})$
for each $t \in [0,1]$,
 \item
$
 (X_1, \omega_1)
$
is isomorphic to $(X, \omega)$
as a K\"{a}hler manifold,
 \item
$\Phi_1$ coincides with $\Phi$ 
under the above isomorphism $X_1 \cong X$,
 \item
$
 (X_0, \omega_0)
$
is a toric variety
with a torus-invariant K\"{a}hler form, 
 \item
$\Phi_0 : X_0 \to \bR^N$ is the moment map
for the torus action on $X_0$, and
 \item 
the flow $\phi_t$ sends $X_{t'}$ to another fiber 
$X_{t'-t}$ preserving the symplectic structures 
and the completely integrable systems:
\[
  \xymatrix{
  (X_{t'}, \omega_{t'})  \ar[dr]_{\Phi_{t'}} 
  \ar[rr]^{ \phi_{t}}
  & & (X_{t'-t}, \omega_{t'-t}) 
      \ar[dl]^{\Phi_{t'-t}} \\ 
  & \bR^N &
            }
          \]
\end{itemize}
\end{df}

Note that
the existence of a toric degeneration
of a projective manifold
with a structure of an integrable system
does not imply the existence of a toric degeneration
of that integrable system.
For example,
$\bP^1 \times \bP^1$ admits
a flat degeneration
into the Hirzebruch surface
$\bF_2 = \bP(\scO_{\bP^1} \oplus \scO_{\bP^1}(2))$,
although the corresponding toric integrable structures
cannot be related by a degeneration
since their moment polytopes are distinct. 

Now the main theorem in this paper is the following:
\begin{thm} \label{th:main}
For any non-increasing sequence
$
 \lambda = (\lambda_1 \ge \lambda_2 \ge \dots \ge \lambda_{n})
$
of real numbers,
the Gelfand-Cetlin system on $(\scO_{\lambda}, \omega_{\lambda})$
admits a toric degeneration.
\end{thm}
\noindent
Essential ingredients of the proof are
the {\em degeneration in stages}
of the flag manifold,
introduced by Kogan and Miller \cite{KM} 
to relate the the Gelfand-Cetlin basis of
$H^0(F, \scL_{\lambda})$
with the monomial basis of $H^0(X_0, \scL_0)$
in a geometric way,
and the {\it gradient-Hamiltonian flow},
introduced by W.-D.~Ruan \cite{R1}
to construct Lagrangian torus fibrations
on Calabi-Yau manifolds.

As an application of Theorem \ref{th:main},
we compute the {\em potential function}
of the Gelfand-Cetlin system
in Theorem \ref{th:potential}
by reducing to the case of toric Fano manifolds,
first studied by Cho and Oh \cite{CO}
and further elaborated by Fukaya, Oh, Ohta and Ono
\cite{FOOO_toric_I}.
The potential function is a Floer theoretic invariant
of a Lagrangian submanifold
introduced by Fukaya, Oh, Ohta and Ono \cite{FOOO2006},
which encodes the information of holomorphic disks
with Lagrangian boundary condition.
It will be used
in Theorem \ref{th:non-displaceable}
to show the existence of a non-displaceable Lagrangian torus
in the flag manifold
just as in the toric case \cite[Theorem 1.5]{FOOO_toric_I}.

In the case of a toric Fano manifold,
the potential function
gives the {\em Landau-Ginzburg potential}
after the substitution of $e^{-1}$
into the indeterminant element $T$ of the Novikov ring
and a suitable change of variables.
The Landau-Ginzburg potential appears
in Givental's integral representation of the {\em $J$-function},
which generates the {\em quantum $D$-module}
encoding the information of Gromov-Witten invariants.
As a corollary to this integral representation,
one obtains an isomorphism between
the quantum cohomology ring
and the Jacobi ring of the Landau-Ginzburg potential.
Such properties continue to hold
in the case of a full flag manifold,
where the $J$-function gives a solution
to the quantum completely integrable system
called the {\em quantum Toda lattice},
although it fails for more general flag manifolds.

The organization of this paper is as follows:
In Section \ref{sc:flag}, we fix notation and
recall basic facts on flag manifolds
which are used through this paper.
In Section \ref{sc:GC-system},
we recall the construction of the Gelfand-Cetlin system.
In Section \ref{sc:toric_deg}, we introduce toric degenerations of flag 
manifolds in stages following \cite{KM}.
A toric degeneration of the Gelfand-Cetlin system is constructed in
Section \ref{sc:deg_of_GC} and Section \ref{Sec:grad-ham}.
In Section \ref{sc:deg_of_GC} we construct a map $\widetilde{\Phi}$ 
in Definition \ref{def:toric_degeneration} using the degeneration
in stages, and prove in Section \ref{Sec:grad-ham} 
that the gradient-Hamiltonian flow sends the flag manifold 
to the toric variety $X_0$ preserving the structure 
of completely integrable systems. 
%
In Section \ref{sc:direct degeneration},
we construct another, not in-stages, toric degeneration
of the Gelfand-Cetlin system
so that $X_t$ is biregular to $X$ for any $t \ne 0$.
This will be used in Section \ref{sc:holo_disk}
to compare the moduli spaces of holomorphic disks
in the flag manifold and the Gelfand-Cetlin toric variety.
In Section \ref{sc:potential},
we recall the definition of the potential function,
and compute it for a Lagrangian torus fiber
in the Gelfand-Cetlin system.
In Section \ref{sc:example},
we study the case of the full flag manifold $F^{(3)}$
and the Grassmannian $\Gr(2, 4)$ in some detail.
In Section \ref{sc:non-displaceable},
we prove the existence of a non-displaceable
Lagrangian torus in the flag manifold
along the lines of \cite{FOOO_toric_I}.
In Section \ref{sc:Toda},
we recall Givental's integral representation
of the $J$-function for the full flag manifold,
and discuss its relation with the potential function.

{\bf Acknowledgment}:
We thank Hiroshi Iritani for
valuable discussions and explanations;
in particular,
we have learned Proposition \ref{prop:valuation} from him.
T.~N. is supported by Grant-in-Aid for Young Scientists (No.19740034)
 and Dean's Grant for Exploratory Research (Graduate School of Science),
 Tohoku University.
Y.~N. is supported by Grant-in-Aid for Young Scientists (No.19740025).
K.~U. is supported by Grant-in-Aid for Young Scientists (No.18840029).

\section{Partial flag manifolds}\label{sc:flag}

Fix a sequence $0 = n_0 < n_1 < \dots < n_r < n_{r+1} =n$ of integers,
and set $k_i = n_i - n_{i-1}$ for $i = 1, \dots , r+1$.
The partial flag manifold $F  = F(n_1 , \dots , n_r , n)$ 
is a complex manifold parameterizing nested subspaces
\[
  0 \subset V_1 \subset \dots \subset V_r \subset \mathbb{C}^n,
  \quad \dim V_i = n_i.
\]
Let $F^{(n)}$ denote the full flag manifold 
$F(1,2, \dots, n)$ for short.
The dimension of $F(n_1, \dots , n_r , n)$ is given by
\[
  N= N(n_1, \dots , n_r, n) := 
  \dim_{\mathbb{C}} F(n_1, \dots , n_r, n) = 
  \sum_{i=1}^r (n_i - n_{i-1})(n - n_i ).
\]
Let $P = P(n_1, \dots , n_r, n) \subset GL(n, \mathbb{C})$ 
be the isotropic subgroup of the standard flag 
$V_i = \langle e_1, \dots , e_{n_i} \rangle$, 
where $\{ e_1, \dots , e_n \}$ is the standard basis of $\mathbb{C}^n$.
Then the intersection of $P$ and $U(n)$ is
$U(k_1) \times \dots \times U(k_{r+1})$, and
$F$ is written as
\[
  F = GL(n,\mathbb{C}) / P
    = U(n) / (U(k_1) \times \dots \times U(k_{r+1}) ).
\]
In particular, the full flag manifold is given by
$F^{(n)} = GL(n,\mathbb{C}) / B = U(n) / T$,
where $B \subset GL(n,\mathbb{C})$ is a Borel subgroup 
consisting of upper triangular invertible matrices, 
and $T$ is a maximal torus in $U(n)$ consisting of 
diagonal matrices.

In this paper we will use two descriptions of flag manifolds, 
(co)adjoint orbits and Pl\"ucker embeddings.
First we recall the (co)adjoint orbit description.
Using a $U(n)$-invariant inner product $\langle \,\, , \, \rangle$ 
on the Lie algebra $\mathfrak u (n)$ of $U(n)$,
we identify the dual $\mathfrak u (n)^*$ of $\mathfrak{u}(n)$ 
with the space $\sqrt{-1} \mathfrak u (n)$ of Hermitian matrices.
We fix $\lambda = \mathrm{diag}\, ( \lambda_1 , \dots , \lambda_n) 
\in \sqrt{-1} \mathfrak u (n)$ with
\begin{equation}
   \underbrace{\lambda_1 = \dots = \lambda_{n_1}}_{k_1}
         > \underbrace{\lambda_{n_1 +1} = \dots = \lambda_{n_2}}_{k_2}
         > \dots > 
        \underbrace{\lambda_{n_r +1} = \dots = \lambda_n}_{k_{r+1}}.
   \label{lambda}
\end{equation}
Then $F$ is identified with the adjoint orbit $\mathcal{O}_{\lambda}
\subset \sqrt{-1} \mathfrak u (n)$ of $\lambda$ by
\[
  F = U(n) / (U(k_1) \times \dots \times U(k_{r+1}) )
  \overset{\sim}{\longrightarrow} \mathcal{O}_{\lambda} ,
  \quad [g] \longmapsto g \lambda g^*.
\]
Note that $\mathcal{O}_{\lambda}$ consists of Hermitian matrices with fixed
eigenvalues $\lambda_1, \dots, \lambda_n$.
$\mathcal{O}_{\lambda}$ has a standard symplectic form $\omega_{\lambda}$
called the Kostant-Kirillov form.
Recall that tangent vectors of $\mathcal{O}_{\lambda}$ at $x$ can be 
written as $\mathrm{ad}_{\xi}(x) = [x, \xi]$ for $\xi \in \mathfrak{u} (n)$.
Then $\omega_{\lambda}$ is defined by
\[
  \omega_{\lambda} 
  \bigl( \mathrm{ad}_{\xi}(x), \mathrm{ad}_{\eta}(x) \bigr) =
  \frac 1{2\pi}
  \langle x , [\xi, \eta] \rangle .
\]
Note that $\omega_{\lambda}$ is the unique $U(n)$-invariant 
K\"ahler form in its cohomology class $[\omega_{\lambda}]$.

Next we recall the Pl\"ucker embedding of $F$.
For each $k = 1 , \dots , n-1$, we set
$\mathbb{P}_k := \mathbb{P} \bigl( \bigwedge^k \mathbb{C}^n \bigr)
 = \mathbb{P}^{\binom nk -1}$.
Then the Pl\"ucker embedding is given by
\[
  \iota : F \hookrightarrow 
  \prod_{i=1}^{r} \mathbb{P}_{n_i},
  \quad 
  (0 \subset V_1 \subset \dots \subset V_r \subset
         \mathbb{C}^n) \mapsto 
  (\textstyle{\bigwedge^{n_1}} V_1, \dots ,
   \textstyle{\bigwedge^{n_r}} V_r).
\]
Note that we have a natural projection
\begin{equation*}
  \begin{matrix}
    \pi = \pi_{n_1,\dots , n_r} : &
    \prod_{k=1}^{n-1} \mathbb{P}_{k} &
    \longrightarrow & \prod_{i=1}^{r} \mathbb{P}_{n_i} \\
    & \cup & & \cup\\
    & F^{(n)} & \longrightarrow & F(n_1 , \dots , n_r,n).
  \end{matrix}
\end{equation*}
For an $n \times n$ matrix $z = (z_{ij})$ and 
$I = \{i_1 < \dots < i_k \} \subset \{1, \dots , n\}$, 
we set
\[
  z_I = \begin{pmatrix} 
        z_{i_1 1} & z_{i_1 2} & \cdots & z_{i_1 k}\\
        z_{i_2 1} & z_{i_2 2} & \cdots & z_{i_2 k}\\
        \vdots    & \vdots    &        & \vdots   \\
        z_{i_k 1} & z_{i_k 2} & \cdots & z_{i_k k}
        \end{pmatrix}.
\]
Then the Pl\"ucker coordinates are given by
\[
  p_I(z) := \det z_I
\]
for $I$ with $|I| = n_1, \dots , n_l$.
In other words, $F$ can be obtained as a ``multiple Proj'' of 
$\mathbb{C}[p_I \, ; \, |I| = n_1, \dots , n_r  \,]$:
\[
  F(n_1, \dots ,n_r, n) = \mathrm{multiple \, Proj}\,
  \mathbb{C}[p_I \, ; \, |I| = n_1, \dots , n_r  \,]
  \subset \prod_{i=1}^{r} \mathbb{P}_{n_i},
\]
which means that $F$ is a subvariety in 
$\prod_{i=1}^{r} \mathbb{P}_{n_i}
 = \prod_{i=1}^{r} \mathrm{Proj}\, \mathbb{C}[Z_I \, ; \, |I| = n_i  \,]$
corresponding to $\mathbb{C}[p_I \, ; \, |I| = n_1, \dots , n_r  \,]$.
In this setting, $\omega_{\lambda}$ coincides with the restriction
$\iota^* \widetilde{\omega}_{\lambda}$ of a K\"ahler form 
\begin{equation}
  \widetilde{\omega}_{\lambda} = \sum_{i=1}^{r}
      (\lambda_{n_i} - \lambda_{n_{i+1}}) 
      \omega_{\mathrm{FS}, n_i}
  \label{Kahler}
\end{equation}
on $\prod_{i=1}^{r} \mathbb{P}_{n_i}$,
where $\omega_{\mathrm{FS}, k}$ is the Fubini-Study form
on $\mathbb{P}_k$.

\begin{expl}
The full flag manifold $F^{(3)}$ for $n=3$ is three dimensional and 
embedded into $\mathbb{P}_1 \times \mathbb{P}_2 = 
\mathbb{P}^2 \times \mathbb{P}^2$ as a hypersurface by
\[
  \iota = 
  ([p_1:p_2:p_3], [p_{12}:p_{13}:p_{23}]) : F^{(3)}
  \longrightarrow \mathbb{P}^2 \times \mathbb{P}^2.
\]
The defining equation (i.e. the Pl\"ucker relation) is given by
\[
  Z_1 Z_{23} - Z_2 Z_{13} + Z_3 Z_{12} = 0,
\]
where $[Z_1:Z_2:Z_3]$, $[Z_{12}:Z_{13}:Z_{23}]$ are homogeneous
coordinates on $\mathbb{P}_1$ and $\mathbb{P}_2$ respectively.
\end{expl}

\begin{expl}
$F^{(4)}$ is of dimension six and embedded into
$\prod_k \mathbb{P}_k = \mathbb{P}^3 \times \mathbb{P}^5
 \times \mathbb{P}^3$.
The Pl\"ucker relations are given by ten quadrics,
and hence $F^{(4)}$ is not a complete intersection.
The projection 
$\pi_2 : \prod_k \mathbb{P}_k \to \mathbb{P}_2 = \mathbb{P}^5$ 
maps $F^{(4)}$ to the Grassmannian $F(2,4) = Gr(2,4)$ 
of two-planes in a four-space, 
which is a hypersurface in $\mathbb{P}_2$ defined by
\[
  Z_{12} Z_{34} - Z_{13} Z_{24} + Z_{14} Z_{23} = 0.
\]
\end{expl}

We consider the case where $\lambda_i \in \mathbb{Z}$ for 
$i = 1, \dots, n$.
Then $\lambda$ can be regarded as a character of $T$ by
\[
  T \longrightarrow \mathbb{C}^*, \quad
  \mathrm{diag}\, (t_1, \dots, t_n) 
  \longmapsto t_1^{\lambda_1} \cdots t_n^{\lambda_n},
\]
and hence gives an action of $T$ on $\mathbb{C}$.
Using this $T$-action, we define a line bundle on $F^{(n)}$ by
\[
  (U(n) \times \mathbb{C})/T \longrightarrow 
  F^{(n)}= U(n)/T.
\]
It is easy to see that this descends to a line bundle $\mathcal{L}_{\lambda}$
on $F(n_1, \dots, n_r,n)$ under the condition (\ref{lambda}).
Note that $\mathcal{L}_{\lambda}$ is also written as
\[
  \mathcal{L}_{\lambda} = \iota^*
     \mathcal{O}_{\mathbb{P}_{n_1}}
     (\lambda_{n_1} - \lambda_{n_2}) 
     \boxtimes \dots \boxtimes 
     \mathcal{O}_{\mathbb{P}_{n_r}}
     (\lambda_{n_r} - \lambda_{n}) ,
\]
where $\boxtimes$ is the outer tensor product.
Hence $\omega_{\lambda}$ represents the first Chern class
$c_1(\mathcal{L}_{\lambda})$ of $\mathcal{L}_{\lambda}$.

We recall the description of the anti-canonical bundle 
$\mathcal{K}_F^{-1}$
on $F$ in terms of characters $\lambda$.
Note that the holomorphic tangent space of $F^{(n)}$ at the 
standard flag can be identified with $\bigoplus_{i<j} 
\mathbb{C} E_{ij}$, where $E_{ij}$ is the matrix whose $(i,j)$-entry 
is 1 and the other entries are zero.
Hence the anti-canonical bundle $\mathcal{K}_{F^{(n)}}^{-1}$ 
of the full flag manifold corresponds to the sum of 
``positive roots'' 
$(0, \dots, 0, 1,0 ,\dots, 0,-1,0,\dots ,0)$:
\[
  \mathcal{K}_{F^{(n)}}^{-1} = L_{2 \rho}, \quad
  2 \rho = (n-1, n-3, \dots, -n+3, -n+1).
\]
Similarly, the anti-canonical bundle of the partial flag manifold
$F(n_1, \dots, n_r,n)$ is given by
\begin{equation}
  \lambda = (\underbrace{n-n_1, \dots}_{k_1}, 
      \underbrace{n-n_1-n_2, \dots}_{k_2}, 
      \dots, 
      \underbrace{n-n_{r-1}-n_r, \dots}_{k_r}, 
      \underbrace{-n_r, \dots, -n_r}_{k_{r+1}} ) ,
  \label{canonical}
\end{equation}
the sum of roots corresponding to $E_{ij}$ above the diagonal squares
$Q_k$ of the ladder diagram, which will be introduced in the following.
For example, the anti-canonical bundle $\mathcal{K}_{Gr(r,n)}^{-1}$ 
of the Grassmannian $Gr(r,n)$ corresponds to
\[
  \lambda = (\underbrace{n-r, \dots, n-r}_{r}, 
      \underbrace{-r, \dots, -r}_{n-r} ).
\]

Now we introduce the standard {\it ladder diagram}, which is used
in \cite{B} to describe the toric degeneration of $F$.
\begin{df}
We consider an $n \times n$ square $Q$ and place squares 
$Q_l$ of size $k_l \times k_l$ ($l=1, \dots, r+1$)
on the diagonal.
The ladder diagram is the set of boxes below the 
diagonal squares.
Let $O_0$ denote the lower left corner of the ladder diagram.
For $l= 1, \dots , r$, the lower right corner of $Q_l$ 
is denoted by $O_l$.
\end{df}
\begin{figure}[h]
 \begin{center}
  \includegraphics*{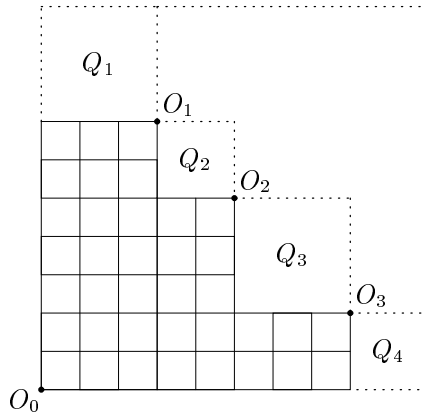}
 \end{center}
 \caption{The ladder diagram for $F(3,5,8,10)$}
 \label{ladder1}
\end{figure}
Note that the number of boxes in the ladder diagram is equal to
$N(n_1, \dots , n_r, n)= \dim_{\mathbb{C}} F(n_1, \dots , n_r, n)$.
Any matrix in $P(n_1, \dots, n_r , n)$ has 0's in its entries which
correspond to boxes in the ladder diagram, and $U(k_l)$ is placed in the 
diagonal square $Q_l$.

\begin{figure}[h]
 \begin{center}
  \includegraphics*{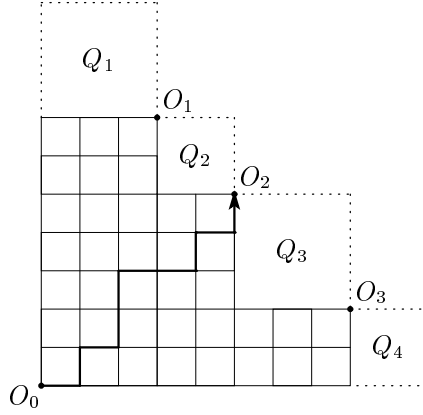}
 \end{center}
 \caption{A positive path}
 \label{ladder2}
\end{figure}

\begin{df}
A {\it positive path} is a path on the ladder diagram, 
starting at the lower left corner $O_0$
and moving either upward or to the right along
edges, until one of $O_k$ is reached.\footnote{%
The choice of the orientation of positive paths is different
from the one in \cite{B}.}
\end{df}
For each positive path ending at $O_k$, we can associate a 
homogeneous coordinate on 
$\mathbb{P}_{n_k} = \mathbb{P}^{\binom n{n_k} -1}$.
Note that the number of positive paths reaching $O_k$ is 
$\binom n{n_k}$.
If the path is horizontal in the $i_1, \dots , i_{n_k}$-th steps,
then the corresponding coordinate is $Z_{i_1, \dots , i_{n_k}}$.
For example, the positive path in Figure \ref{ladder2} corresponds
to $Z_{1,3,6,7,9}$ on $\mathbb{P}_5$.

\section{The Gelfand-Cetlin system}\label{sc:GC-system}

In this section we recall the construction of the Gelfand-Cetlin system.
First we consider the case of full flag manifold 
$F^{(n)} = F(1,2,\dots , n)$.
Note that this corresponds to the case where all $\lambda_i$ 
are distinct:
\[
  \lambda_1 > \lambda_2 > \dots > \lambda_n.
\] 
For $x \in \mathcal{O}_{\lambda}$ and $k = 1, \dots ,n-1$, 
let $x^{(k)}$ denote the upper-left $k \times k$ submatrix of $x$.
Since $x^{(k)}$ is also a Hermitian matrix, it has real eigenvalues
$\lambda^{(k)}_1(x) \ge \lambda^{(k)}_2(x) \ge \dots 
\ge \lambda^{(k)}_k(x)$.
By taking the eigenvalues for all $k = 1, \dots ,n-1$,
we obtain a set of functions
\begin{equation}
  \Phi_{\lambda} : \mathcal{O}_{\lambda} \longrightarrow 
  \mathbb{R}^{n(n-1)/2}, \quad
  x \longmapsto \left( \lambda^{(k)}_i (x) \right)_{%
    \begin{subarray}{l}
    k = 1, \dots , n-1,\\
    i = 1, \dots ,k
    \end{subarray}}.
  \label{G-Csystem}
\end{equation}
Recall that $\dim_{\mathbb{C}} F^{(n)} = n(n-1)/2$.

\begin{thm}[Guillemin-Sternberg \cite{GS}]\label{Th:G-S}
$\{\lambda_i^{(k)} \}_{k,i}$ is a completely
integrable system on $(\mathcal{O}_{\lambda}, \omega_{\lambda})$. 
\end{thm}

$\{\lambda_i^{(k)} \}_{k,i}$ is called the 
{\it Gelfand-Cetlin system} on $(F^{(n)} , \omega_{\lambda})$.

\begin{rem}
For $k=1, \dots , n-1$, we regard $U(k)$ as a subgroup of $U(n)$ by
\[
  U(k) \cong \left( \begin{array}{c|c}
             U(k) & 0 \\
            \hline  0 &1_{n-k}
          \end{array} \right)
  \subset U(n).
\]
Then the map $x \mapsto x^{(k)}$ gives a moment map of 
the $U(k)$-action on $(\mathcal{O}_{\lambda}, \omega_{\lambda})$.
\end{rem}

For later use, we present a proof of this theorem.
We first recall some basic facts on moment maps.
Let $G$ be a compact Lie group acting on a symplectic manifold
$(M, \omega)$ with a moment map $\mu : M \to \mathfrak{g}^*$.
Note that $\mathfrak{g}^*$ has a Poisson structure induced from the 
Kostant-Kirillov form.
For functions $f_1$, $f_2$ on $\mathfrak{g}^*$,
their Poisson bracket $\{ f_1 , f_2 \}_{\mathfrak{g}^*}$ at 
$x \in \mathfrak{g}^*$ is defined to be the Poisson bracket at $x$ 
of the restrictions $f_i|_{\mathcal{O}_x}$ to the coadjoint
orbit $\mathcal{O}_x$ of $x$.
\begin{lem}\label{Lem:collective}
For $f_1, f_2 \in C^{\infty}(\mathfrak{g}^*)$, it follows that
\[
  \{ \mu^* f_1, \mu^* f_2 \}_{M}
  = \mu^* \{ f_1, f_2 \}_{\mathfrak{g}^*},
\]
where $\{ \,\, , \, \}_M$ is the Poisson bracket on $M$.
In particular, if $f_1$ \rm{(}or $f_2$\rm{)} is 
$\mathrm{Ad}(G)^*$-invariant, then we have
\[
  \{ \mu^* f_1, \mu^* f_2 \}_{M} = 0.
\]
\end{lem}

See \cite{GS2} for a proof.
We also recall the following Noether type theorem, which will 
be used in Section \ref{Sec:grad-ham}.
\begin{lem}\label{Lem:Noether}
  If $f \in C^{\infty}(M)$ is $G$-invariant, then 
  $\mu$ is constant along the Hamiltonian flow of $f$.
\end{lem}

Now we go back to our situation and prove Theorem \ref{Th:G-S}.
Since $\lambda_i^{(k)}$ is a pull back of a $U(k)$-invariant function
on $\sqrt{-1} \mathfrak u (k)$ by the moment map
$x \mapsto x^{(k)}$,
Lemma \ref{Lem:collective} implies that
\[
  \{ \lambda_i^{(k)}, \lambda_j^{(l)} \} =0
\]
on $(\mathcal{O}_{\lambda}, \omega_{\lambda})$.

Next we see the image $\Phi_{\lambda} (\mathcal{O}_{\lambda})$ 
of $\mathcal{O}_{\lambda}$.
We first consider the eigenvalues of $x$ and $x^{(n-1)}$.
The mini-max principle implies that
\begin{equation*}
  \lambda_1 \ge \lambda^{(n-1)}_1 \ge \lambda_2 \ge 
  \lambda^{(n-1)}_2 \ge \lambda_3 \ge \dots \ge 
  \lambda_{n-1} \ge \lambda^{(n-1)}_{n-1} \ge \lambda_n.
\end{equation*}
Hence $\left( \lambda^{(k)}_i (x) \right)$ satisfies
\begin{equation}
\begin{alignedat}{17}
  \lmd 1 &&&& \lmd 2 &&&& \lmd 3 && \cdots && \lmd {n-1} &&&& \lmd n  \\
  & \uge && \dge && \uge && \dge &&&&&& \uge && \dge & \\
  && \llmd {n-1}1 &&&& \llmd {n-1}2 &&&&&&&& \llmd{n-1}{n-1} && \\
  &&& \uge && \dge &&&&&&&& \dge &&& \\
  &&&& \llmd {n-2}1 &&&&&&&& \llmd{n-2}{n-2} &&&& \\
  &&&&& \uge &&&&&& \dge &&&&& \\
  &&&&&& \dndots &&&& \updots &&&&&& \\
  &&&&&&& \uge && \dge &&&&&&& \\
  &&&&&&&& \llmd 11 &&&&&&&&& 
\end{alignedat}
\label{GC-pattern}
\end{equation}
An array of real numbers satisfying (\ref{GC-pattern}) is called
a {\it Gelfand-Cetlin pattern} for $\lambda$.
The {\it Gelfand-Cetlin polytope} $\Delta_{\lambda}$ is a polytope
consisting of Gelfand-Cetlin patterns for $\lambda$.
The above argument means the image  
$\Phi_{\lambda} (\mathcal{O}_{\lambda})$ is
contained in $\Delta_{\lambda}$.

\begin{lem}
Let $a_1, \dots , a_{k+1}, b_1, \dots, b_k$ be real numbers 
satisfying
\[
  a_1 \ge b_1 \ge a_2 \ge  \dots \ge a_k \ge b_k \ge a_{k+1}.
\]
Then there exist $x_1, \dots, x_k \in \mathbb{C}$ and 
$x_{k+1} \in \mathbb{R}$ such that
\[
  \begin{pmatrix} b_1 & & 0 & \bar{x}_1\\
                  & \ddots & & \vdots\\
                  0 & & b_k & \bar{x}_k\\
                  x_1 & \hdots & x_k & x_{k+1}
  \end{pmatrix}
\]
has eigenvalues $a_1, \dots , a_{k+1}$.
\end{lem}

We omit the proof.
Using this lemma successively, we can prove that 
$\Phi_{\lambda} (\mathcal{O}_{\lambda})= \Delta_{\lambda}$.
The fact that $\dim \Delta_{\lambda} = n(n-1)/2$ implies the
functional independence of $\lambda_i^{(k)}$'s.

\begin{rem}
  Recall that a moment map of the $T$-action on 
  $\mathcal{O}_{\lambda}$ is given by
  \[
    x = (x_{ij}) \longmapsto \begin{pmatrix}
                               x_{11} & & 0\\
                               & \ddots & \\
                               0 & & x_{nn}
                             \end{pmatrix}.
  \]
  Since 
  \[
    x_{kk} = \tr x^{(k)} - \tr \, x^{(k-1)}
         = \sum_i \lambda^{(k)}_i - \sum_i \lambda^{(k-1)}_i, 
  \] 
  the $T$-action is contained in the Gelfand-Cetlin system.
\end{rem}

\begin{rem}
We can think of the ladder diagram as a container of a Gelfand-Cetlin 
pattern leaning to the right, here 
$\lambda_1, \dots , \lambda_n$ are placed in the diagonal 
squares $Q_1, \dots, Q_n$ respectively
(see Figure \ref{ladder3}).
\begin{figure}[h]
  \begin{center}
    \includegraphics*{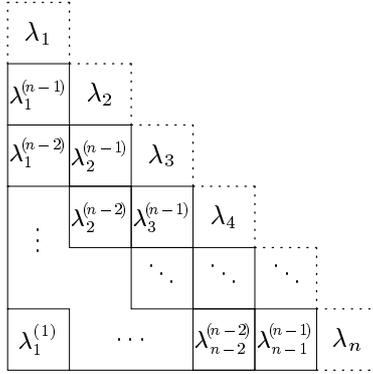}
  \end{center}
  \caption{The ladder diagram as a container of Gelfand-Cetlin patterns}
  \label{ladder3}
\end{figure}
A vertex of the Gelfand-Cetlin polytope is given by 
a Gelfand-Cetlin pattern each of whose entry is 
connected to some $\lambda_i$ by a chain of equalities.
By putting arrows on edges of the ladder diagram where 
the adjacent entries are distinct,
we obtain a tree of positive paths, which is called a
{\it meander} in \cite{B}.
\end{rem}

\begin{expl} \label{eg:F3}
In the case of $F^{(3)}$, the Gelfand-Cetlin system consists of 
three functions
\[
  \Phi_{\lambda} = (\lambda_1^{(2)}, \lambda_2^{(2)}, \lambda_1^{(1)})
  : F^{(3)} \longrightarrow 
  \mathbb{R}^3,
\]
and the Gelfand-Cetlin polytope is illustrated in Figure \ref{GCpoly}.
\begin{figure}[h]
 \begin{center}
  \includegraphics*{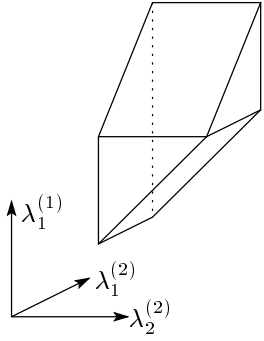}
 \end{center}
 \caption{The Gelfand-Cetlin polytope for $F^{(3)}$}
 \label{GCpoly}
\end{figure}
Topology of the fibers are quite similar to the toric case: 
for almost every point in $\Delta_{\lambda}$, 
if it is contained in an $i$-dimensional face, 
its fiber is an $i$-dimensional torus.
A difference appears at the vertex where four edges are 
intersecting. 
The fiber of this point is a three dimensional sphere $S^3$.
This can be seen as follows.
The vertex is given by the equations
\[
  \begin{matrix}
    \lambda_1 &&&& 
    \hspace{-10pt} \lambda_2 &&&& 
    \hspace{-10pt} \lambda_3 \\
    & \phantom{\hspace{-10pt} \lueq} && \hspace{-10pt} \rueq 
    && \hspace{-10pt} \lueq  && 
      \phantom{\hspace{-10pt} \rueq}&\\
    && \hspace{-10pt} \hbox to1.05em{$\lambda_1^{(2)}$} &&&&
    \hspace{-10pt} \hbox to1.05em{$\lambda_2^{(2)}$}  && \\
    &&& \hspace{-10pt} \lueq && \hspace{-10pt} \rueq &&& \\
  &&&& \hspace{-10pt} \hbox to1.05em{$\lambda_1^{(1)}$} &&&& 
  \end{matrix}.
\]
Hence a point $x \in \mathcal O_ {\lambda}$ in the fiber of this point
must satisfy
\[
  x^{(2)} = \begin{pmatrix}
             \lambda_2 & 0 \\
             0 & \lambda_2
            \end{pmatrix},
\]
or equivalently, $x$ must have the form
\[
  x = \begin{pmatrix}
        \lambda_2 &  & z_1 \\
         & \lambda_2 & z_2 \\
        \bar{z}_1 & \bar{z}_2 & \nu
      \end{pmatrix}.
\]
The condition that $x$ has eigenvalues $\lambda_1 , \lambda_2, \lambda_3$
is equivalent to 
\[
  \nu = \lambda_1 - \lambda_2 + \lambda_3, \quad
  |z_1|^2 + |z_2|^2 = (\lambda_1 - \lambda_2)
                      (\lambda_2 - \lambda_3),
\]
which means that the fiber is isomorphic to $S^3$.
\end{expl}

\begin{rem}
The torus action given by the Gelfand-Cetlin system does not preserve
the complex structure on $F$.
In fact, any torus acting holomorphically on $F$ must contained in a 
maximal torus of $U(n)$.
Thus the inverse image of a face of the Gelfand-Cetlin polytope is not
necessarily a complex subvariety in $F$.
In the case of $F^{(3)}$, for two faces of dimension two 
in the back side of $\Delta_{\lambda}$ in Figure \ref{GCpoly}, 
their inverse images are complex subvarieties.
On the other hand, it is not true for other four faces of dimension two.
\end{rem}

We move on to the case of a partial flag manifold 
$F(n_1, \dots, n_r, n) \cong \mathcal{O}_{\lambda}$ 
where $\lambda$ satisfies (\ref{lambda}).
We can consider the functions (\ref{G-Csystem}) also in this case.
Under the condition (\ref{lambda}), (\ref{GC-pattern}) implies that
\begin{align*}
  \lambda_1^{(n-1)} = \dots = &\lambda_{n_1-1}^{(n-1)}
  = \lambda_{n_1},\\
  \lambda_{n_1+1}^{(n-1)} = \dots = &\lambda_{n_2-1}^{(n-1)}
  = \lambda_{n_2},\\
  &\vdots 
\end{align*}
which mean that $\lambda_i^{(k)}$ contained in $Q_l$ is a constant
function $\lambda_{n_l}$.
In other words, non-constant $\lambda_i^{(k)}$ exactly corresponds to
a box in the ladder diagram.
In particular, we have the right number of Poisson commuting functions
\[
  \Phi_{\lambda} : (F(n_1, \dots, n_r,n), \omega_{\lambda}) 
  \longrightarrow 
  \mathbb{R}^{N(n_1, \dots, n_r,n)}, \quad
  x \longmapsto \left( \lambda^{(i)}_j (x) \right).
\]
We call this the Gelfand-Cetlin system on $F(n_1, \dots, n_r, n)$.
\begin{figure}[h]
 \begin{center}
  \includegraphics*{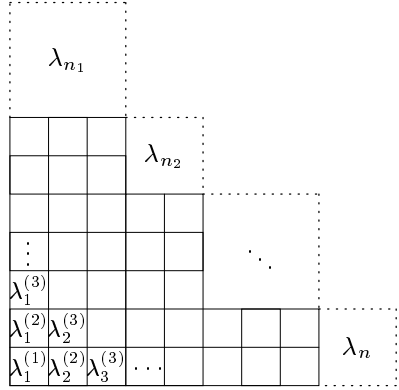}
 \end{center}
 \caption{Gelfand-Cetlin patterns in a partial flag case}
 \label{ladder4}
\end{figure}

\begin{expl}
We consider the case of $Gr(2,4)$, where the condition
$\lambda_1 = \lambda_2 > \lambda_3 = \lambda_4$ is satisfied.
The Gelfand-Cetlin system consists of four functions
\[
  \Phi_{\lambda} = (\lambda_2^{(3)}, \lambda_1^{(2)}, 
    \lambda_2^{(2)}, \lambda_1^{(1)})
  : Gr(2,4) \longrightarrow \mathbb{R}^4.
\]
Under the projection $\mathbb{R}^4 \to \mathbb{R}$, 
$(\lambda_2^{(3)}, \lambda_1^{(2)}, \lambda_2^{(2)}, \lambda_1^{(1)})
\to \lambda_2^{(3)}$, 
$\Delta_{\lambda}$ is fibered by Gelfand-Cetlin polytopes for $F^{(3)}$,
and the fiber shrinks to a two dimensional triangle
on the boundaries $\lambda_2^{(3)} = \lambda_1, \lambda_3$.
We see a fiber of a boundary point given by
$\lambda_2^{(3)} = \lambda_1^{(2)} = \lambda_2^{(2)} = \lambda_1^{(1)}
 = s$.
From the argument in the $F^{(3)}$-case, 
each matrix $x$ in this fiber satisfies
\[
  x^{(3)} = \begin{pmatrix}
        s &  & z_1 \\
         & s & z_2 \\
        \bar{z}_1 & \bar{z}_2 & \alpha
      \end{pmatrix}
\]
with $\alpha = \lambda_1 + \lambda_3 - s$ and 
$|z_1|^2 + |z_2|^2 = (\lambda_1 - s)(s - \lambda_3)$.
We assume that $\lambda_1 > s > \lambda_3$, and
take $g \in U(3)$ such that $g^* x^{(3)} g 
= \mathrm{diag}(\lambda_1, s, \lambda_3)$.
Note that such $g$ is unique up to scalar.
Then it is easy to check that $x$ has the form
\[
  \begin{pmatrix} g & \\ & 1 \end{pmatrix}^* x
  \begin{pmatrix} g & \\ & 1 \end{pmatrix}
  = \begin{pmatrix} \lambda_1 &&& 0\\
                              & s && z_3\\
                              && \lambda_3 & 0\\
                              0 & \bar{z}_3 & 0 & \alpha \end{pmatrix}
\] 
with $|z_3|^2 =(\lambda_1 - s)(s - \lambda_3)$.
This means that the fiber is a Lagrangian $S^3 \times S^1$.
When $s$ goes to a boundary, say $\lambda_1$, we have
$x = \mathrm{diag}( \lambda_1, \lambda_1, \lambda_3, \lambda_3)$,
which means that the fiber shrinks to a point.
\end{expl}

In the rest of this section, we see some properties of the Gelfand-Cetlin polytope.
Readers who are interested only in toric degeneration of the Gelfand-Cetlin system
can skip to the next section.
First we consider the case where $\omega_{\lambda}$ represents
the anti-canonical class, or equivalently, $\lambda$ is given by
(\ref{canonical}).

\begin{df}\label{def:reflexive}
An $N$-dimensional integral polytope $\Delta \subset \mathbb{R}^N$
is said to be {\it reflexive} if the following two conditions fold:
\begin{enumerate}
\renewcommand{\labelenumi}{$(\mathrm{\roman{enumi}})$}
 \item all codimension one faces of $\Delta$ are supported by an
       affine hyperplane of the form
       $\{ u \in \mathbb{R}^N \, | \, \langle u,v \rangle = -1 \}$
       for some $v \in \mathbb{Z}^N$, where $\langle \, , \, \rangle$
       is the standard inner product on $\mathbb{R}^N$.
 \item $\Delta$ contains only one integral point $0$ in its interior.
\end{enumerate}
\end{df}

It is proved in \cite{B1} that $\Delta$ is reflexive if and only if 
the corresponding toric variety is Fano.

\begin{lem}
Under the condition $(\ref{canonical})$, the Gelfand-Cetlin polytope
$\Delta_{\lambda}$ is reflexive after a translation.
\end{lem}

See also \cite{B}, where the same is proved in the dual side.

\begin{proof}
For each $k$ and $i$, we set $\lambda^{(k)}_i = k - 2i + 1$.
Then it is easy to check that the collection of these 
$\lambda_i^{(k)}$ gives a Gelfand-Cetlin pattern.
The definition implies that
\[
  \lambda_i^{(k+1)} = \lambda_i^{(k)} +1 > \lambda_i^{(k)}
  > \lambda_i^{(k)} -1 = \lambda_{i+1}^{(k+1)} 
\]
for each $k$ and $i$,
which means that this $(\lambda_i^{(k)})$ is the unique 
integral point in the interior of $\Delta_{\lambda}$.
\begin{figure}[h]
 \begin{center}
  \includegraphics*{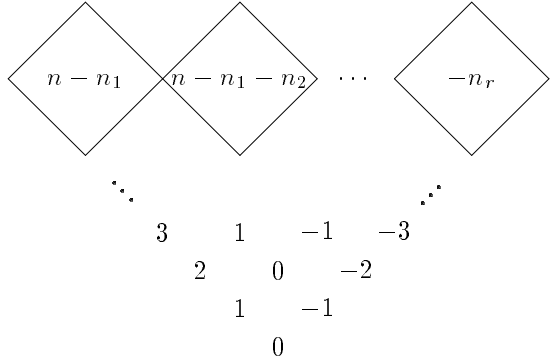}
 \end{center}
 \caption{An integral point in the Gelfand-Cetlin polytope}
 \label{ladder7}
\end{figure}
\end{proof}

Next we compute the volume of the Gelfand-Cetlin polytope $\Delta_{\lambda}$, 
though it is not used in the following.

\begin{prop}
Under the condition $(\ref{lambda})$, the volume of the
Gelfand-Cetlin polytope $\Delta_{\lambda}$ is given by
\[
  \mathrm{Vol}(\Delta_{\lambda})
  = \frac{ \prod_{i<j} (\lambda_{n_i} - \lambda_{n_j})^{k_i k_j}}
    {\prod_{k=1}^{n-1} k! }.
\]
\end{prop}

\begin{proof}
When all $\lambda_i$ are integers, the Borel-Weil theory states that
the space $H^0(F,\mathcal{L}_{\lambda})$ of holomorphic sections of 
$\mathcal{L}_{\lambda}$
is an irreducible representation of $U(n)$ of highest weight $\lambda$.
Gelfand and Cetlin \cite{GC} constructed a basis of 
$H^0(F,\mathcal{L}_{\lambda})$
called the Gelfand-Cetlin basis, which is indexed by integral
points of $\Delta_{\lambda}$.
In particular, the dimension of $H^0(F,\mathcal{L}_{\lambda})$ is equal to the
number $\# (\Delta_{\lambda} \cap \mathbb{Z}^N)$ of integral points in
$\Delta_{\lambda}$.
On the other hand, the dimension of $H^0(F,\mathcal{L}_{\lambda})$ is given by
the Weyl dimension formula
\[
  \dim H^0(F,\mathcal{L}_{\lambda}) 
  = \frac{ \prod_{i<j} (\lambda_i - \lambda_j + j -i)}
    {\prod_{k=1}^{n-1} k! }.
\]
The proposition follows from
\[
  \frac 1{m^N} \# (\Delta_{m \lambda} \cap \mathbb{Z}^N)
  = \mathrm{Vol}(\Delta_{\lambda}) + O \left( \frac 1m \right)
\]
and 
\begin{align*}
  \frac 1{m^N} \dim H^0(F,\mathcal{L}_{m \lambda}) 
  &= \frac 1{m^N} 
    \frac{ \prod_{i<j} (m \lambda_i - m \lambda_j + j -i)}
    {\prod_{k=1}^{n-1} k! }\\
  &= \frac{ \prod_{i<j} (\lambda_{n_i} - \lambda_{n_j})^{k_i k_j}}
    {\prod_{k=1}^{n-1} k! }
    + O \left( \frac 1m \right)
\end{align*}
for sufficiently large $m \in \mathbb{Z}_{>0}$.
\end{proof}

For example, in the full flag case, the volume of $\Delta_{\lambda}$ is given
be the difference product of $\lambda_1, \dots, \lambda_n$:
\[
  \mathrm{Vol}(\Delta_{\lambda})
  = \frac{ \prod_{i<j} (\lambda_i - \lambda_j)}{\prod_{k=1}^{n-1} k! }.
\]

We close this section with computation of the volume of the dual polytope 
$(\Delta_{\lambda})^*$ of $\Delta_{\lambda}$ for full flag manifolds and 
Grassmannians, which will be used in Section \ref{sc:non-displaceable}.
First we consider the case of full flag manifolds.
Let $e^{(k)}_j$ denote the unit vector corresponding to the coordinate $\lambda^{(k)}_j$.
Then $(\Delta_{\lambda})^*$ is a convex hull of $\pm e^{(n-1)}_j$, $e^{(k+1)}_j - e^{(k)}_j$, and
$e^{(k)}_j - e^{(k+1)}_{j+1}$.

\begin{lem}
  In the full flag case, the volume of $(\Delta_{\lambda})^*$ is equal to  $2^N/N!$, where 
  $N= n(n-1)/2 = \dim F^{(n)}$.
\end{lem}

To see this, we observe that $\Delta^*$ can be constructed successively as follows.
We start with $\pm e^{(n-1)}_1$, which give an interval $[-1,1]$.
Adding $\pm e^{(n-1)}_2$ we obtain a square which is a union of two triangles 
of hight 1 over the base $[-1,1]$, glued along the interval.
By adding $\pm e^{(n-1)}_3$ further, we get an octahedron, 
two pyramids of hight 1 over the square,
glued along the common base square, and so on.
The construction is almost the same for rows below $\lambda^{(n-1)}_j$:
adding $e^{(k+1)}_j - e^{(k)}_j$ and $e^{(k)}_j - e^{(k+1)}_{j+1}$
we also obtain two cones of hight 1 over the polytope constructed at this stage, 
glued along their bases.
Therefore the volume of $(\Delta_{\lambda})^*$ is given by
\[
  2 \cdot \frac 22 \cdot \frac 23 \cdots \frac 2N = \frac {2^N}{N!}.
\]

\begin{figure}[h]
 \begin{center}
  \includegraphics*{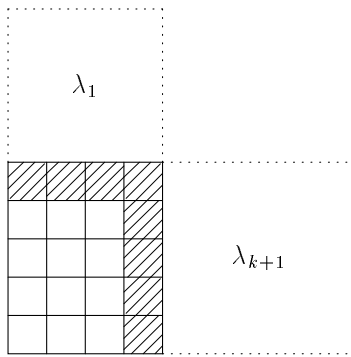}
 \end{center}
 \caption{A Gelfand-Cetlin pattern for a Grassmannian}
 \label{ladder9}
\end{figure}
In the case of Grassmannian $Gr(k,n)$, we start with the shaded boxes 
in Figure \ref{ladder9}. 
Since 
\[
  \lambda_1 \ge u_1 \ge u_2 \ge \dots \ge u_{n-1} \ge \lambda_{k+1}
\]
are the only nontrivial relations, this part gives a simplex with vertices
\[
  \begin{pmatrix} -1 \\ 0 \\ 0 \\ \vdots \\ \\  0 \end{pmatrix},
  \begin{pmatrix}  1 \\ -1 \\ 0 \\ \vdots \\ \\  0 \end{pmatrix},
  \begin{pmatrix}  0 \\ 1 \\ -1 \\ 0 \\ \vdots \\  0 \end{pmatrix},
  \dots,
  \begin{pmatrix}  0 \\ 0 \\ \vdots \\ 0 \\ 1 \\  -1 \end{pmatrix},
  \begin{pmatrix}  0 \\ 0 \\ \vdots \\  \\ 0 \\ 1 \end{pmatrix}.
\]
Its volume is given by
\[
  \frac 1{(n-1)!} \det \begin{pmatrix}
      2 &  1 &         1 & \hdots &  1 & 1\\     
     -1 &  1 &         0 &           &      &   \\
      0 & -1 &         1 &         0 &      &   \\
         &  \ddots   & \ddots & \ddots &  \ddots   &   \\
         &     &         0  &       -1 &   1 &  0 \\
       0 &     &            &        0 &  -1 & 1
    \end{pmatrix}
  = \frac n{(n-1)!}.
\]
Since contributions of the remaining boxes are the same as in 
the full flag case, we obtain the following.
\begin{lem}
  The volume of $(\Delta_{\lambda})^*$ for the Grassmannian $Gr(k,n)$ is given by
  \[
    \mathrm{Vol}(\Delta_{\lambda})^* = \frac {n 2^{N-(n-1)}}{N!},
  \]
  where $N = k(n-k) = \dim Gr(k,n)$.
\end{lem}

\section{Refinement of Gelfand-Cetlin polytopes and small resolutions}
\label{sc:resolution}

In this section, we show the following result
about the geometry of Gelfand-Cetlin polytopes:
\begin{prop}\label{refinement}
Every refinement of the fan $\Sigma$ of a Gelfand-Cetlin toric variety
 into simplicial cones, without adding
 new rays, gives a resolution of the Gelfand-Cetlin toric variety.
In particular, this gives a small resolution.
\end{prop} 
\begin{rem}
This statement was shown in \cite{Batyrev}
for a particular refinement of the fan.
See also the discussion in the last part of Section \ref{sc:toric_deg} 
for the full flag case.
Below we give a different proof, which gives the result for any refinement.
\end{rem}

To prove the proposition, it suffices to show the following lemma.
\begin{lem}\label{lem:monoid}
   Let $\sigma$ be any $N$-dimensional simplicial cone given by 
    N-rays of a maximal dimensional cone in the fan $\Sigma$.
   Then the monoid $\bZ^N \cap \sigma$ is generated by the 
    integral generators of the rays of $\sigma$.
\end{lem}
\proof
Since each ray in $\Sigma$ corresponds to a toric 
 divisor, which is given by Gelfand-Cetlin patterns
 such that just one of the inequalities is chosen to be an
 equality, the choice of rays can be expressed 
 as a set of equalities in a Gelfand-Cetlin pattern.
Each equality has the form 
\[
 \lambda_i - \lambda_i^{(n-1)} = 0
 \quad \text{or} \quad
 \lambda_i^{(n-1)} - \lambda_{i+1} = 0
\] 
 if the equality appears in the top row 
 of the Gelfand-Cetlin pattern, and 
 \[
   \lambda_i^{(k-1)} - \lambda_{i+1}^{(k)} = 0
   \quad \text{or} \quad
   \lambda_i^{(k)} - \lambda_i^{(k-1)} = 0
 \] 
 otherwise.
Hence the generators of one-dimensional cones in $\Sigma$ have the form
 $(0, \dots, 0, \pm 1, 0, \dots, 0)$ or 
 $(0, \dots, 0, \pm 1, 0, \dots, 0, \mp 1,0, \dots, 0)$.
Recall that each cone in $\Sigma$ of maximal dimension corresponds to a 
 vertex of the Gelfand-Cetlin polytope, and the vertex is given by 
 a Gelfand-Cetlin pattern such that every entry is connected to some 
 $\lambda_i$ in the top row by a chain of equalities.
Note also that every singular locus is a toric stratum such that the
 corresponding Gelfand-Cetlin patterns contain a loop of equalities 
 such as
 \[
   \begin{alignedat}{5}
    && \llmd{k+1}{i+1} && \\
    & \deq && \ueq & \\
    \llmd{k}{i} &&&& \llmd{k}{i+1} \\
    & \ueq && \deq & \\
    && \llmd{k-1}{i} &&,
   \end{alignedat}
 \]
 and such a loop makes the corresponding cone to be 
 non-simplicial.
Hence choosing $N$-rays which give a simplicial cone is equivalent to
 removing some of the equalities in the Gelfand-Cetlin pattern
 so that the resulting set of equalities does not form any loop.
Note that the resulting chain of equalities may contain a part
which does not occur in Gelfand-Cetlin patterns such as
\[
   \begin{alignedat}{5}
    && \llmd{k+1}{i+1} && \\
    & \deq && \ueq & \\
    \llmd{k}{i} &&&& \llmd{k}{i+1} \\
    & \ueq && \dge & \\
    && \llmd{k-1}{i} &&.
   \end{alignedat}
 \]
By arraying the generators of the rays of $\sigma$
 considered as column vectors,
 we obtain an $N \times N$ matrix $A$.
Then Lemma \ref{lem:monoid} follows from the following:
\begin{claim}
  The matrix $A$ has determinant $\pm 1$.
\end{claim}
\proof
Since the chain of equalities above is a tree when it is regarded 
 as a graph in such a way that its edges are given by the equalities,
 we can take a univalent end.
The row in $A$ corresponding to this univalent end
 has the form $(0, \dots, 0, \pm1, 0, \dots, 0)$,
 and hence the calculation of the determinant can be reduced 
 to that for an $(N-1) \times (N-1)$ matrix.
We obtain the above claim by repeating this process.
\qed\\
      
\section{Degeneration of flag manifolds in stages} \label{sc:toric_deg}

It is known that $F(n_1, \dots , n_r, n)$ degenerates into the 
{\it Gelfand-Cetlin toric variety}, 
a toric variety which corresponds
to the Gelfand-Cetlin polytope \cite{GL, KM, B}.
In this section, we recall the construction of toric degenerations
in \cite{KM}, with minor changes.

The toric degeneration is given by deforming the Pl\"ucker embedding.
For that purpose, 
we introduce a weight $w_{ij}$ of each variable $z_{ij}$ given by
\[
  w_{ij} = \begin{cases}
      3^{i-j-1}, & i > j,\\
      0,           & i \le j.
    \end{cases}
\]
Namely, the matrix of weights $w_{ij}$ is given by
\[
  w = (w_{ij}) = \begin{pmatrix}
             0       &        &        &        &  \\
             1       & 0      &        &        &  \\
             3       & 1      & 0      &        &  \\
             \vdots  & \ddots & \ddots & \ddots &  \\
             3^{n-2} & \ldots & 3      & 1      & 0
           \end{pmatrix}.
\]
For each $I = \{ i_1 < \dots < i_k \} \subset \{1, \dots, n \}$, 
we set
\begin{equation*}
  q_I(z,t) := t^{- \tr w_I} p_I( t^{w_{ij}}z_{ij})
  = t^{- \tr \, w_I} \det ( t^{w_{ij}}z_{ij})_I.
\end{equation*}
Since the diagonal term
\[
  d_I (z) = z_{i_1 1} z_{i_2 2} \dots z_{i_k k}
\]
of $p_I( z_{ij}) = \det z_I$ is the unique term of the lowest weight, 
$q_I(z,t)$ is a polynomial in $z_{ij}$ and $t$.
From the construction, $q_I(z,1) = p_I(z)$ is a Pl\"ucker 
coordinate for $t=1$, and 
$q_I(z,0) = d_I (z)$ is a {\it monomial} for $t=0$.
We define a one-parameter family of projective varieties by
\[
  \mathfrak{X} = \mathfrak{X}(n_1, \dots, n_l,n)
  = \mathrm{multiple \, Proj} \, \mathbb{C} 
  [t, q_I \, ; \, |I| = n_1, \dots , n_r  \,]
  \subset
  \prod_{i=1}^r \mathbb{P}_{n_i} \times \mathbb{C}.
\]
In the full flag case, we simply write
$\mathfrak{X}(1, \dots,n) = \mathfrak{X}^{(n)}$.
Note that $\mathfrak{X}(n_1, \dots,n_r,n)$ is obtained from 
$\mathfrak{X}^{(n)}$ by the natural projection
\[
  \begin{array}{cccl}
  \pi_{n_1, \dots, n_l} : &
  \mathfrak{X}^{(n)} & \longrightarrow
  & \mathfrak{X}(n_1, \dots,n_r,n)\\
  & \downarrow & & \downarrow \\
  & \mathbb{C} & = & \mathbb{C} .
  \end{array}
\]
\begin{thm}
  $f : \mathfrak{X} \to \mathbb{C}$ is a flat family of projective 
  varieties such that $X_1 := f^{-1}(1)$ is the flag manifold $F$, 
  and the central fiber $X_0 := f^{-1}(0)$ is the 
  Gelfand-Cetlin toric variety.
\end{thm}

The existence of toric degenerations is first proved by Gonciulea 
and Lakshmibai \cite{GL},
and the fact that $X_0$ is isomorphic to the Gelfand-Cetlin toric variety
is proved by Kogan and Miller \cite{KM} in the full flag case.
The results are generalized to partial flag cases by Batyrev et al. \cite{B0, B}.
Note that $X_0$ is singular except for the trivial case, i.e. the case of 
projective spaces.
It is proved in \cite{B} that the singular locus of $X_0$ consists of 
codimension three conifold strata.

\begin{rem}
  For every $t \ne 0$, $X_t = f^{-1}(t)$ is isomorphic to the
  flag manifold $F$ as a complex manifold.
  On the other hand, the restriction 
  $\widetilde{\omega}_{\lambda} |_{X_t}$ of the K\"ahler form
  defined in (\ref{Kahler}) coincides with
  the Kostant-Kirillov form only when $|t|=1$.
  Note also that the natural $U(k)$-action on $\prod \mathbb{P}_{n_i}$ ($k \ge 2$)
  does not preserves $X_t$ for $|t| \ne 1$ in general.
\end{rem}

\begin{expl}
The degenerating family for the full flag manifold
$F^{(3)}$ of dimension three is given by 
\[
  \mathfrak{X} =
  \Bigl\{ 
    \bigl([Z_1 : Z_2 : Z_3], [Z_{12} : Z_{13} : Z_{23}] ,t \bigr)
   \, \Bigm| Z_1 Z_{23} - Z_2 Z_{13} + t Z_3 Z_{12} =0 \,
   \Bigr\},
\]
with the central fiber
\[
  X_0 = \Bigl\{ 
    \bigl([Z_1 : Z_2 : Z_3], [Z_{12} : Z_{13} : Z_{23}] \bigr)
    \in \mathbb{P}^2 \times \mathbb{P}^2 
   \, \Bigm|  Z_1 Z_{23} - Z_2 Z_{13}  = 0 \,
   \Bigr\}.
\]
Note that $X_0$ has a singularity at $([0:0:1], [1:0:0])$,
which corresponds to the vertex of the Gelfand-Cetlin polytope
emanating four edges.
\end{expl}

\begin{expl}
Recall that $Gr(2,4)$ is embedded into $\mathbb{P}^5$ with its defining 
equation $Z_{12} Z_{34} - Z_{13} Z_{24} + Z_{14} Z_{23} = 0$.
The equation for the toric degeneration of $Gr(2,4)$ is given by
\[
  t Z_{12} Z_{34} - Z_{13} Z_{24} + Z_{14} Z_{23} = 0.
\]
Then $X_0$ has conifold singularities along 
$\{ Z_{13} = Z_{24} = Z_{14} = Z_{23} = 0 \} = \mathbb{P}^1$.
\end{expl}

\begin{figure}[h]
 \begin{center}
  \includegraphics*{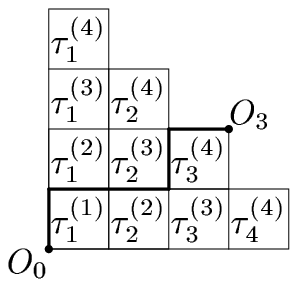}
  \caption{Positive paths and monomials on $X_0$.}
  \label{ladder5}
 \end{center}
\end{figure}

Now we see the Gelfand-Cetlin toric variety $X_0$ more closely.
Let $I = \{ i_1 < \dots < i_k \}$ and 
$J = \{ j_1 < \dots < j_l \}$ be subsets of $\{1, \dots, n\}$ with 
$k \le l$, and $\gamma_I$, $\gamma_J$ the corresponding positive 
paths.
We define their meet and join by
\begin{align*}
  I \wedge J &= \{ \min (i_1, j_1), \dots , \min (i_k, j_k),
                   j_{k+1}, \dots , j_l \},\\
  I \vee J&= \{\max (i_1, j_1), \dots , \max (i_k, j_k) \}.
\end{align*}
Then $I \wedge J$ (resp. $I \vee J$) corresponds to a positive path
moving along the lower (resp. upper) route of the union 
$\gamma_I \cup \gamma_J$.
The defining equations for the Gelfand-Cetlin toric variety
$X_0 \subset \prod_i \mathbb{P}_{n_i}$ are
given by the following binomial relations
\begin{equation}
  Z_I Z_J - Z_{I \wedge J} Z_{I \vee J} = 0
  \label{binomial}
\end{equation}
(see \cite{GL} or \cite{KM}).
Next we see the monomial embedding of $X_0$ into 
$\prod \mathbb{P}_{n_i}$.
Let $T_k$ be a torus corresponding to the $k$-th row  
$( \lambda_i^{(k)} )_i$ from the bottom of Gelfand-Cetlin patterns.
In the full flag case, $T_k$ is a $k$-dimensional torus $T^k$.
We take natural coordinates $( \tau_i^{(k)})_i$ on 
$T_k^{\mathbb{C}} = T_k \otimes \mathbb{C}$,
and consider the following matrix
\[
  \tau = \begin{pmatrix}\begin{array}{llllll}
      \tau^{(n-1)}_1 \dots \tau^{(2)}_{1} \tau^{(1)}_{1} &
      &&&&\\
      \tau^{(n-1)}_1 \dots \tau^{(2)}_{1} &
      \tau^{(n-1)}_2 \dots \tau^{(2)}_{2}
      &&&&\\
      \vdots & \vdots & \ddots &\\
      \tau^{(n-1)}_1 \tau^{(n-2)}_1 &
      \tau^{(n-1)}_2 \tau^{(n-2)}_2
      & \ldots & \tau^{(n-1)}_{n-2} \tau^{(n-2)}_{n-2} &&\\
      \tau^{(n-1)}_1 & \tau^{(n-1)}_2
      & \ldots & \tau^{(n-1)}_{n-2} & \tau^{(n-1)}_{n-1} &\\
      1 & 1 & \ldots & 1 & 1 & 1
    \end{array}\end{pmatrix},
\]
where we assume that $\tau^{(k)}_i = 1$ if the corresponding
$\lambda_i^{(k)}$ is contained in a diagonal square $Q_l$.
Then the embedding of $X_0$ into 
$\prod_i \mathbb{P}_{n_i}$ is given by the monomials
\begin{equation}
  Z_{I} = d_I (\tau), \quad
  |I|= n_1, \dots, n_r.
\label{monomial}
\end{equation}
These monomials can be described using the ladder diagram in the following way.
We put $\tau^{(k)}_i$ on the ladder diagram in the same way 
as for $\lambda_i^{(k)}$.
Then each monomial is written as
$d_I(\tau) = \prod \tau_i^{(k)}$, where the product is taken
over $\tau_i^{(k)}$'s placed above the positive path $\gamma_I$
corresponding to $I$.
For example, the path in Figure \ref{ladder5} corresponds to
$\tau^{(4)}_1 \tau^{(3)}_1 \tau^{(2)}_1 \tau^{(4)}_2 \tau^{(3)}_2$.
It is easy from this expression to see that these monomials satisfy the
binomial relations (\ref{binomial}).

We extend $f: \mathfrak{X} \to \mathbb{C}$ to an $(n-1)$-parameter 
family to define the {\it degeneration of flag manifolds in stages}
which is also introduced in \cite{KM},
on which we will construct a toric degeneration of the Gelfand-Cetlin system.
Let $\boldsymbol{t} = (t_2, \dots , t_n)$ be parameters 
and set
\[
  \widetilde{w}_{k,ij} = \begin{cases}
      0, & i < k,\\
      w_{kj} - w_{k-1,j},  & i \ge k
    \end{cases}
\]
for $k = 2, \dots, n$.
Then the weight of $z_{ij}$ is extended to a multi-weight
\[
  \boldsymbol{t}^{\widetilde{w}_{ij}} :=
  t_2^{\widetilde{w}_{2,ij}} t_3^{\widetilde{w}_{3,ij}}
  \cdots t_n^{\widetilde{w}_{n,ij}}.
\]
Then the matrix of multi-weights is given by
\[
 ( \boldsymbol{t}^{\widetilde{w}_{ij}} )_{ij}
  = \begin{pmatrix}\begin{array}{llllll}
      1 &&&&&\\
      t_2 & 1 &&&&\\
      t_2 t_3^2 & t_3 & 1 &\\
      t_2 t_3^2 t_4^6 & t_3 t_4^2 & t_4 & 1 \\
      \vdots & \vdots &  & \ddots & \ddots &\\
      t_2 t_3^2 t_4^6 \dots t_n^{2 \cdot 3^{n-3}} &
      t_3 t_4^2 \dots t_n^{2 \cdot 3^{n-4}} & \ldots &
      &  t_n & 1
    \end{array}\end{pmatrix}.
\]
Note that $t_k$ does not appear above the $k$-th row.
We set
\[
  \tilde{q}_I (z_{ij},t_2, \dots , t_n) = 
  d_I (\boldsymbol{t}^{\widetilde{w}_{ij}})^{-1}
  p_I( \boldsymbol{t}^{\widetilde{w}_{ij}} z_{ij})
  \in \mathbb{C} [z_{ij}, t_k ].
\]
Since $\boldsymbol{t}^{\widetilde{w}_{ij}} = t^{w_{ij}}$ 
for $\boldsymbol{t} = (t, \dots, t)$, we have
$\tilde{q}_I (z_{ij},t, \dots ,t) =q_J(z_{ij}, t)$.
By taking multiple Proj of $\mathbb{C} [t_i, \tilde{q}_I]$, 
we obtain an $(n-1)$-parameter family 
$\tilde{f} : \widetilde{\mathfrak{X}} \to \mathbb{C}^{n-1}$
of projective varieties.
From the construction, 
$X_{(1,\dots,1)} := \tilde{f}^{-1}(1,\dots,1)$ is isomorphic 
to the flag manifold $F$, and
$X_{(0,\dots,0)} = \tilde{f}^{-1}(0,\dots,0)$ is the 
Gelfand-Cetlin toric variety $X_0$.
An important point is that the $U(k-1)$-action on
$\prod_i \mathbb{P}_{n_i}$ preserves each fiber 
$X_{(1, \dots 1, t_k, \dots, t_n)}$ for 
$t_2 = \dots = t_{k-1} = 1$,
and the action of $T_{n-1} \times \dots \times T_{k}$ preserves
$X_{(t_2, \dots, t_k, 0, \dots , 0)}$ for 
$t_{k+1} = \dots = t_{n} = 0$
(see \cite{KM}, or the discussion below).
Hence we consider a sequence of degenerations given 
by varying the parameters as follows:
\[
  \boldsymbol{t} = 
  (1, \dots , 1) \rightsquigarrow (1, \dots,1, 0) \rightsquigarrow 
      (1, \dots, 1, 0, 0) \rightsquigarrow \dots \rightsquigarrow 
      (0, \dots , 0). 
\]
Let
\begin{equation}
  \begin{matrix}
    f_k : & \mathfrak{X}_k = 
    \widetilde{\mathfrak{X}}|_{\begin{subarray}{l}
     t_2 = \dots = t_{k-1} = 1\\
     t_{k+1} = \dots = t_n =0
     \end{subarray}} 
    & \longrightarrow & \mathbb{C}\\
    & \cup & & \rotatebox[origin=c]{90}{$\in$} \\
    & X_{(1, \dots , 1, t_k , 0, \dots , 0)} & \longrightarrow & t_k 
  \end{matrix}
  \label{stage}
\end{equation}
denote the $(n-k+1)$-th stage of the degeneration, i.e. a one-parameter 
sub-family given by fixing $t_2 = \dots = t_{k-1} = 1$ 
and $t_{k+1} = \dots = t_n =0$.
We write $X_{k, t_k} := 
f_k^{-1}(t_k) = X_{(1, \dots , 1, t_k , 0, \dots , 0)}$ for short.
Note that each $X_{k,t}$ has actions of  
$T_{n-1} \times \dots \times T_k$ and $U(k-1)$.
In particular, $X_{k,1} = X_{k+1,0}$ admits actions of
$T_{n-1} \times \dots \times T_{k}$ and $U(k)$.

\begin{rem}
 The final stage $f_2 : \mathfrak{X}_2 \to \mathbb{C}$  of the 
 toric degeneration is a trivial family.
 For example, the two-parameter family for $F^{(3)}$ is given by
 \[
    t_2 Z_1 Z_{23} - t_2 Z_2 Z_{13} + t_3^2 Z_3 Z_{12} = 0.
 \]
 This fact is related to the fact 
 that $F^{(2)} = \mathbb{P}^1$ is a toric variety and the family 
 in this case is trivial (see the discussion below).
\end{rem}

Now we consider the full flag case, and see 
$X_{n, 0} = X_{n-1,1} = X_{(1,\dots,1,0)}$ in more detail,
where the actions of $T_{n-1}$ and $U(n-1)$ can be seen 
explicitly as follows:
\begin{lem}
$X_{n-1,1}$ has a resolution 
$h_{n-1,1}:Y_{n-1,1} \to X_{n-1,1}$ such that
$Y_{n-1,1}$ has a structure of $(\mathbb{P}^1)^{n-1}$-bundle
over a smaller flag manifold $F^{(n-1)}$.
The $U(n-1)$-action on $X_{n-1,1}$ is induced from the standard
one on $F^{(n-1)}$, and the $T_{n-1}$-action comes from the 
natural torus action on the $(\mathbb{P}^1)^{n-1}$-fibers.
\end{lem}

\begin{proof}
We consider the multi-parameter family
$\widetilde{\mathfrak{X}}^{(n-1)} \to \mathbb{C}^{n-2}$ 
for $F^{(n-1)}$.
Recall that each fiber is a  subvariety in 
$\mathbb{P}^{(n-1)} := 
  \prod_k \mathbb{P}(\textstyle{\bigwedge^k} \mathbb{C}^{n-1})$.
Let $\mathcal{E}_0 = \mathcal{E}_{n-1} 
= \mathcal{O}_{\mathbb{P}^{(n-1)}}$ be trivial bundles and
$\mathcal{E}_k = \mathrm{pr}_k^* 
\mathcal{O}_{\mathbb{P}(\bigwedge^k \mathbb{C}^{n-1})}(1)$
for each $k$, where
$
  \mathrm{pr}_k : \mathbb{P}^{(n-1)} 
  \to \mathbb{P}(\textstyle{\bigwedge^k} \mathbb{C}^{n-1})
$
is the natural projection.
We define a $(\mathbb{P}^1)^{n-1}$-bundle on $\mathbb{P}^{(n-1)}$ by 
\[
  E :=
  \mathbb{P} (\mathcal{E}_0 \oplus \mathcal{E}_1) \times_{\mathbb{P}^{(n-1)}}
  \mathbb{P} (\mathcal{E}_1 \oplus \mathcal{E}_2) \times_{\mathbb{P}^{(n-1)}}
  \dots \times_{\mathbb{P}^{(n-1)}}
  \mathbb{P} (\mathcal{E}_{n-2} \oplus \mathcal{E}_{n-1}).
\]
Restricting this to each fiber $X_{(t_2, \dots, t_{n-1})}^{(n-1)}$,
we obtain a family
\[
  \begin{matrix}
  \mathfrak{Y} &\longrightarrow &\widetilde{\mathfrak{X}}^{(n-1)}\\
  \cup & & \cup \\
  Y_{(t_2, \dots, t_{n-1})} & \longrightarrow 
   & X_{(t_2, \dots, t_{n-1})}^{(n-1)}
  \end{matrix}
\]
of $(\mathbb{P}^1)^{n-1}$-bundles.
We claim that there exists a surjective birational morphism
\[
  \begin{matrix}
  \mathfrak{Y} &\longrightarrow &\widetilde{\mathfrak{X}}|_{t_n=0}\\
  \cup & & \cup \\
  Y_{(t_2, \dots, t_{n-1})} & \longrightarrow & X_{(t_2, \dots, t_{n-1},0)} \,.
  \end{matrix}
\]
To see this, we observe that
\begin{equation}
  \tilde{q}_I(z, t_2, \dots, t_{n-1},0) =
  \begin{cases} 
  \tilde{q}_I(z^{(n-1)}, t_2, \dots, t_{n-1}),
  & i_k < n,\\
  z_{nk} \tilde{q}_{i_1 , \dots , i_{k-1}}(z^{(n-1)}, 
  t_2, \dots, t_{n-1}),
  & i_k = n
  \end{cases}
  \label{1st-deg}
\end{equation}
for $I= \{ i_1 < \dots < i_k\} \subset \{1, \dots, n\}$,
where $\tilde{q}_I$ in the right hand side is regarded as a function
of $(n-1) \times (n-1)$ matrices $z^{(n-1)}$. 
Let $Z'_I$, $I \subset \{1, \dots, n-1 \}$ be the homogeneous 
coordinates of $\prod_k \mathbb{P} (\bigwedge^k \mathbb{C}^{n-1})$,
where we assume that $Z'_{\emptyset} = Z'_{1, \dots, n-1} = 1$, 
and $[u_i:v_i]$ 
the fiber coordinates of $\mathbb{P} (\mathcal{E}_{i-1} \oplus \mathcal{E}_i)$
with $u_i \in \mathcal{E}_{i-1}$ and $v_i \in \mathcal{E}_i$.
Then (\ref{1st-deg}) implies that
\begin{equation}
  Z_I = \begin{cases}
    u_k Z'_I, \quad \text{if $n \not\in I$},\\
    v_k Z'_{i_1, \dots, i_{k-1}}, \quad \text{if $I = \{ i_1, \dots, i_{k-1}, 
    i_k = n \}$}
  \end{cases}
  \label{birat}
\end{equation}
on $\mathbb{P}_k$ gives a surjective birational morphism
\[
  Y_{(t_2, \dots , t_{n-1})} \subset E
  \longrightarrow 
  X_{(t_2, \dots , t_{n-1},0)} \subset \prod \mathbb{P}_k
\]
for each $(t_2, \dots , t_{n-1})$
(we may think that $z_{nk} = v_k/u_k$).
In particular, we have a resolution 
\[
  h_{n-1,1}: Y_{n-1,1} := Y_{(1, \dots , 1)}
  \longrightarrow 
  X_{n-2, 1} = X_{(1, \dots , 1,0)}
\]
of $X_{n-1,1}$.
We can see from (\ref{birat}) that the $U(n-1)$-action on $X_{n-1, 1}$ 
is induced from the standard one on 
$F^{(n-1)} = X^{(n-1)}_{(1, \dots , 1)}$.
On the other side, (\ref{monomial}) implies that the $T_{n-1}$-action
on $\prod \mathbb{P}_k$ is given by
\[
  Z_I \longmapsto  \begin{cases}
    \tau_1^{(n-1)} \dots \tau_k^{(n-1)} Z_I, \quad 
      \text{if $n \not\in I$},\\
    \tau_1^{(n-1)} \dots \tau_{k-1}^{(n-1)} Z_I, \quad 
    \text{if $I = \{ i_1, \dots, i_{k-1}, i_k = n \}$}
  \end{cases}
\]
on $\mathbb{P}_k$, or equivalently,
\[
  [Z_I: Z_{I'n}]_{I, I' \subset \{1, \dots, n-1\}} 
  \longmapsto 
  [\tau_k^{(n-1)} Z_I: Z_{I'n}]_{I, I' \subset \{1, \dots, n-1\}}.
\]
This together with (\ref{birat}) mean that the action of $T_{n-1}$ 
is induced from a natural action on the fiber $(\mathbb{P}^1)^{n-1}$.
\end{proof}

By repeating this process, 
we have an iterated fibration
\[
  Y_{m,1} \overset{(\mathbb{P}^1)^{n-1}}{\longrightarrow}
  \dots \overset{(\mathbb{P}^1)^m}{\longrightarrow}
  F^{(m)}
\]
over $F^{(m)}$ and a resolution 
$h_{m,1} : Y_{m,1} \longrightarrow X_{m,1} $
for each $m = n-1, \dots , 2$
such that the $U(m)$-action on $X_{m,1}$ is induced from 
the standard one on $F^{(m)}$, and the torus action comes from 
a natural action on the fibers.
Taking the Pl\"ucker coordinates $([Z'_I]_{|I|=k})_k$ on $F^{(m)}$ and
fiber coordinates 
$([u^n_i : v^n_i])_{i=1, \dots, n-1} \in (\mathbb{P}^1)^{n-1}, 
\dots,
([u^{m+1}_i : v^{m+1}_i])_{i=1,\dots, m+1} \in (\mathbb{P}^1)^m$, 
$h_{m,1}$ is given by
\begin{equation}
  Z_I = \chi_{l, I''}(u,v) \cdot Z'_{I'}
  \label{resol}
\end{equation}
for $I= \{ i_1 < \dots < i_k < i_{k+1} < \dots < i_l \}
\subset \{1, \dots , n\}$ with $I'= \{ i_1 , \dots ,i_k \}
\subset \{ 1, \dots , m \}$ and 
$I'' =\{i_{k+1} , \dots , i_l \} \subset \{ m+1, \dots , n\}$,
where $\chi_{l,I''}(u,v)$ is a monomial in $u^k_i$, $v^l_j$ defined
inductively by (\ref{birat}).
$\chi_{l,I''}(u,v)$ is explicitly given by
\[
  \chi_{l, I''}(u,v) =
  u^n_l \dots u^{i_l+1}_l v^{i_l}_l u^{i_l-1}_{l-1}
  \dots u^{i_{l-1}+1}_{l-1} v^{i_{l-1}}_{l-1}
  u^{i_{l-1}-1}_{l-2} \dots v^{i_{k+1}}_{k+1}
  u^{i_{k+1}-1}_k \dots u^{m+1}_k .
\]
Note that (\ref{resol}) also gives a birational surjective morphism
$h_{m,t} : Y_{m,t} \to X_{m,t}$ for each $t$.

In particular, we obtain a resolution $Y_0 \to X_0$
of the Gelfand-Cetlin toric variety such that $Y_0$ has a
structure of iterated fibration over $\bP^1$.
This is a small resolution of $X_0$ constructed in
\cite{B}.

\section{Toric degeneration of Gelfand-Cetlin systems}\label{sc:deg_of_GC}

In this and the next sections we construct a toric degeneration of 
the Gelfand-Cetlin system using the degeneration in stages 
discussed in the previous section.

For each $m = 1, \dots, n-1$, we consider the natural $U(m)$-action on 
$\prod_{j=1}^r \mathbb{P}_{n_j}$ which is an extension of
the action on $F(n_1, \dots, n_r,n)$.
Let
\[
  \mu^{(m)} : \prod_{j=1}^r \mathbb{P}_{n_j} \longrightarrow
  \sqrt{-1} \mathfrak{u} (m)
\]
denote the moment map of the $U(m)$-action, and
\[
  \tilde{\lambda}_i^{(m)} : \prod_{j=1}^r \mathbb{P}_{n_j} \longrightarrow
  \mathbb{R}, \quad
  \tilde{\lambda}_1^{(m)} \ge \dots \ge \tilde{\lambda}_m^{(m)}
\]
be functions which associate eigenvalues of $\mu^{(m)}(Z)$ 
for each $Z  \in \prod_{j=1}^r \mathbb{P}_{n_j}$.
From the construction, the collection of
$\tilde{\lambda}_j^{(m)}$'s restricted to $X_1$ gives the Gelfand-Cetlin 
system on $F(n_1, \dots, n_r,n)$.
Hereafter we do not assume that indices of $Z_I$ are increasing, 
and use the convention
\[
  Z_{\sigma (i_1) , \dots, \sigma (i_k)} = 
  (\mathrm{sgn} \, \sigma ) Z_{i_1, \dots, i_k}
\]
for $\sigma \in \mathfrak{S}_k$.
Then the moment map $\mu^{(n)}$ of the $U(n)$-action is given by
\[
  \mu^{(n)} (Z) = \sum_{k=1}^{r} 
  \frac{\lambda_{n_k} - \lambda_{n_{k+1}}}{\sum_{|I|=n_k} |Z_I|^2}
  \left( \sum_{|I'|=n_k-1} Z_{i I'} \overline{Z}_{j I'} 
  \right)_{i,j=1, \dots, n}
  + \lambda_n \cdot 1_n,
\]
and $\mu^{(m)}$ is its $m \times m$ upper left block
\begin{equation}
  \mu^{(m)} (Z) = \sum_{k=1}^{r} 
  \frac{\lambda_{n_k} - \lambda_{n_{k+1}}}{\sum_{|I|=n_k} |Z_I|^2}
  \left( \sum_{|I'|=n_k-1} Z_{i I'} \overline{Z}_{j I'} 
  \right)_{i,j= 1, \dots, m}
  + \lambda_n \cdot 1_m.
  \label{mu}
\end{equation}
We also extend to $\prod_{j=1}^r \mathbb{P}_{n_j}$ the torus action 
on the Gelfand-Cetlin toric variety $X_0$,
and consider the moment map
\[
  \tilde{\nu}^{(m)}_i : \prod_{j=1}^r \mathbb{P}_{n_j} 
  \longrightarrow \mathbb{R}
\]
of the action of $\tau^{(m)}_i$.
From (\ref{monomial}), the torus action is given by
\[
  Z_I \longmapsto d_I(\tau) \cdot Z_I.
\]
In particular, the action of $T_m$ is given by
\[
  Z_I \longmapsto \tau_1^{(m)} \dots \tau_k^{(m)} Z_I
\]
for $I = \{ i_1 < \dots < i_k < i_{k+1} < \dots < i_l \}$ with
$\{ i_1 < \dots < i_k\} \subset \{1, \dots , m\}$ and
$\{ i_{k+1} < \dots < i_l \} \subset \{ m+1, \dots, n\}$,
where we assume that $\tau^{(m)}_i = 1$ if it is contained in a
diagonal square $Q_j$.
Hence we have
\begin{equation}
  \tilde{\nu}^{(m)}_j = \sum_{k=1}^r 
  \frac {\lambda_{n_k} - \lambda_{n_{k+1}}}{\sum_{|I|=n_k} |Z_I|^2}
  \sum_{\begin{subarray}{c}
    |I|=n_k, \\ i_j \le m \end{subarray}} 
  |Z_I|^2
  + \lambda_n .
  \label{nu}
\end{equation}
Note that 
$\tilde{\lambda}_1^{(1)} = \tilde{\nu}^{(1)}_1$.
For each $t$, we define
\[
  \Phi_{m,t} = \Bigl.
  \left( \tilde{\nu}_i^{(n-1)}, \dots, \tilde{\nu}_j^{(m)},
  \tilde{\lambda}_k^{(m-1)}, \dots, \tilde{\lambda}_1^{(1)}
  \right) \Bigr|_{X_{m,t}}
  : X_{m,t}
  \longrightarrow \mathbb{R}^{N(n_1, \dots, n_r,n)}.
\]
Then $\Phi_{n,1}$ coincides with the Gelfand-Cetlin system
on $X_{n,1} = F$, while $\Phi_{2,0}$ is the moment map 
of the torus action on the Gelfand-Cetlin toric variety $X_0$.

\begin{thm}\label{Thm:CIS}
Let $f_m: \mathfrak{X}_m \to \mathbb{C}$ be the $(n-m+1)$-th stage
of the toric degeneration defined in $(\ref{stage})$ 
$(m =n, \dots , 2)$.
Then, for each $t$, $\Phi_{m,t} : X_{m,t} \to \mathbb{R}^{N}$
is a completely integrable system on
$(X_{m,t}, \widetilde{\omega}_{\lambda}|_{X_{m,t}})$.
Moreover, for $t = 0$,
\[
  \tilde{\lambda}_k^{(m-1)}
   =
  \tilde{\nu}_k^{(m-1)}
\]
on $X_{m,0} = X_{m-1,1}$.
Namely, $\Phi_{m,0}$ coincides with the initial map $\Phi_{m-1,1}$
in the next stage.
\end{thm}

\begin{proof}
Functional independence for $\tilde{\lambda}^{(k)}_i$'s and
$\tilde{\nu}^{(l)}_j$'s follows from the fact that, for all $t$, the image 
$\Phi_{m,t}(X_{m,t})$ coincides with the Gelfand-Cetlin polytope 
$\Delta_{\lambda}$, whose dimension is equal to 
$\frac 12 \dim_{\mathbb{R}} X_{m,t}$.
This will be proved in the next section
(Corollary \ref{Cor:flow}).
We prove here the Poisson commutativity of the functions.
Since the symplectic structure and the moment maps 
$\mu^{(k)}$, $\tilde{\nu}^{(l)}_j$ in the full flag case descend to 
those in the partial flag case under the condition (\ref{lambda}),
it suffices to prove in the full flag case.

Since the $U(m-1)$-action preserves $X_{m,t}$, 
the restriction $\mu^{(m-1)}|_{X_{m,t}}$ gives a moment map
of the $U(m-1)$-action on $X_{m,t}$.
From Lemma \ref{Lem:collective}, 
we have
\[
  \{ \tilde{\lambda}_i^{(k)}, \tilde{\lambda}_j^{(l)} \} =0,
  \quad k, l \le m-1
\]
on $X_{m,t}$.
Similarly, the restrictions of $\tilde{\nu}^{(k)}_i$ to $X_{m,t}$, 
$(k \ge m)$ give a moment map
of the $T_{n-1}\times \dots \times T_m$-action on $X_{m,t}$,
and hence we obtain
\[
  \{ \tilde{\nu}^{(k)}_i, \tilde{\nu}^{(l)}_j \} = 0
\]
for $k,l \ge m$ on $X_{m,t}$.

To see that $\tilde{\lambda}_i^{(k)}$ commutes with 
$\tilde{\nu}_j^{(l)}$ ($k < m$, $l \ge m$), we pull them back
to $Y_{m,t}$.
Hereafter we normalize homogeneous coordinates so that
$\sum_{|I|=k} |Z_I|^2 = \sum_{|I|=k} |Z'_I|^2 = 
 |u_i^k|^2 + |v_i^k|^2 = 1$,
and assume that $\lambda_n = 0$ for simplicity.
From (\ref{mu}) and (\ref{resol}) we have
\begin{align*}
  h_m^* \mu^{(m)} 
  &= \sum_{l=1}^{n-1} (\lambda_l - \lambda_{l+1})
     \Biggl( \sum_{k=1}^l \sum_{%
       \begin{subarray}{c}
         I' \subset \{1, \dots, m\},\\
         |I'| = k-1,\\
         I'' \subset \{ m+1, \dots, n\},\\
         |I''| = l-k
       \end{subarray}}
     \chi_{l, I''} Z'_{iI'} \overline{\chi_{l, I''} Z'_{jI'}}
    \Biggr)_{i,j=1, \dots, m}\\
  &= \sum_{l=1}^{n-1} (\lambda_l - \lambda_{l+1})
     \left\{ \sum_{k=1}^l \Biggl( \sum_{%
       \begin{subarray}{c} 
         I'' \subset \{ m+1, \dots, n\},\\
         |I''| = l-k
       \end{subarray}}
       |\chi_{l,I''}|^2 \Biggr)
     \Biggl( \sum_{%
       \begin{subarray}{c} 
         I' \subset \{1, \dots, m\},\\
         |I'| = k-1
       \end{subarray}}
       Z'_{iI'} \overline{Z}'_{jI'} \Biggr)_{i,j=1, \dots, m}
     \right\}.
\end{align*}
Since $T_l$ acts only on the fibers $([u^j_i: v^j_i])_i$, 
$h_m^* \mu^{(m)}$ (and hence $h_m^* \mu^{(k)}$, $k < m$) is 
$T_l$-invariant.
In particular, we have
\[
  \{ \tilde{\lambda}_i^{(k)}, \tilde{\nu}_j^{(l)} \} 
  = \xi_{\tilde{\nu}_j^{(l)}} \tilde{\lambda}_i^{(k)} =0,
  \quad k \le m-1, \quad l \ge m
\]
on $X_{m,t}$,
where $\xi_{\tilde{\nu}_j^{(l)}}$ is the Hamiltonian vector field
of $\tilde{\nu}_j^{(l)}$.

Finally we check that $\tilde{\lambda}_j^{(m)} 
= \tilde{\nu}_j^{(m)}$ on $X_{m,1} = X_{m+1,0}$,
which is also seen on the resolution $Y_{m,1}$.
Since $h_m^* \mu^{(m)}$ has a form of a moment map 
of the standard $U(m)$-action on $F^{(m)}$,
the eigenvalues $h_{m,1}^* \tilde{\lambda}_j^{(m)}$ are $U(m)$-invariant.
In particular, 
$h_{m,1}^* \tilde{\lambda}_j^{(m)}$ is determined 
by its values on the fiber over the standard flag in $F^{(m)}$.
Recall that the standard flag is given by
\[
  Z'_I = \begin{cases}
         1, \quad \text{if $I = \{1, \dots, k\}$,}\\
         0, \quad \text{otherwise}
        \end{cases}
\]
on $\mathbb{P}_k$.
Then we have
\[
  h_m^* \mu^{(m)} = 
  \sum_{l=1}^{n-1} (\lambda_l - \lambda_{l+1})
     \Bigg\{ \sum_{k=1}^l \Biggl( \sum_{%
       \begin{subarray}{c} 
         I'' \subset \{ m+1, \dots, n\},\\
         |I''| = l-k
       \end{subarray}}
       |\chi_{l,I''}|^2 \Biggr)
     \begin{pmatrix}
            1_k & \\
             & 0_{m-k}
      \end{pmatrix} \Bigg\}
\] 
on the fiber of this point.
Thus its eigenvalues are given by
\begin{equation}
  h_m^* \tilde{\lambda}^{(m)}_j = \sum_{l=1}^{n-1}
  (\lambda_l - \lambda_{l+1})
     \Bigg\{ \sum_{k=j}^l \Biggl( \sum_{%
       \begin{subarray}{c} 
         I'' \subset \{ m+1, \dots, n\},\\
         |I''| = l-k
       \end{subarray}}
       |\chi_{l,I''}|^2 \Biggr)
     \Bigg\}.
  \label{lambda2}
\end{equation}
On the other hand, from (\ref{nu}) we have
\begin{align*}
  h_m^* \tilde{\nu}^{(m)}_j 
  &= \sum_{l=1}^{n-1} (\lambda_l - \lambda_{l+1})
    \Biggl( \sum_{k=j}^l \sum_{%
       \begin{subarray}{c}
         I' \subset \{1, \dots, m\},\\
         |I'| = k,\\
         I'' \subset \{ m+1, \dots, n\},\\
         |I''| = l-k
       \end{subarray}} | \chi_{l, I''} Z'_{I'}|^2
    \Biggr)\\
  &= \sum_{l=1}^{n-1} (\lambda_l - \lambda_{l+1})
    \Biggl\{ \sum_{k=j}^l \Bigg( \sum_{%
       \begin{subarray}{c}
         I'' \subset \{ m+1, \dots, n\},\\
         |I''| = l-k
       \end{subarray}} | \chi_{l, I''} |^2 \Bigg)
       \Bigg( \sum_{%
       \begin{subarray}{c}
         I' \subset \{1, \dots, m\},\\
         |I'| = k,
       \end{subarray}} | Z'_{I'}|^2 \Bigg)
    \Biggr\}
\end{align*}
Since $ \sum_{|I'| = k} |Z'_{I'}|^2 =1$, it follows that
\[
  h_m^* \tilde{\nu}^{(m)}_j 
  = \sum_{l=1}^{n-1} (\lambda_l - \lambda_{l+1})
    \Biggl( \sum_{k=j}^l \sum_{%
       \begin{subarray}{c}
         I'' \subset \{ m+1, \dots, n\},\\
         |I''| = l-k
       \end{subarray}} | \chi_{l, I''}|^2
    \Biggr),
\]
which coincides with (\ref{lambda2}).
\end{proof}

\section{Gradient-Hamiltonian flows}\label{Sec:grad-ham}

For a family of hypersurfaces in a K\"ahler manifold,
W.-D. Ruan \cite{R1} constructed a flow, 
called the {\it gradient-Hamiltonian flow},
which sends a member of the family to another.
We apply this to the toric degenerations of flag manifolds in this section
and finish the construction of a toric degeneration
of the Gelfand-Cetlin system from the previous section.

Let $(\mathfrak X,\widetilde{\omega})$ be a K\"ahler manifold.
Assume that we have a family $X_t = \{ f=t \}$ $(t \in \mathbb{C})$ 
of complex hypersurfaces defined by a meromorphic function 
$f$ on $\mathfrak X$.
For example, if $X_t$ is given by 
$X_t = \{ s_0 - t s_{\infty} = 0 \}$ for some holomorphic sections 
$s_0, s_{\infty} \in H^0(\mathfrak X ,\mathfrak L)$ of 
a line bundle $\mathfrak L$ on $\mathfrak X$, 
then we choose $f= s_0/s_{\infty}$.
Let $\nabla (\re f)$ be the gradient vector field of the
real part of $f$, and
$\xi_{\im f}$ the Hamiltonian vector field of the imaginary part
of $f$.
Then the Cauchy-Riemann equation implies that
\begin{equation*}
  \nabla (\re f) = - \xi_{\im f}.
\end{equation*}
We define the {\it gradient-Hamiltonian vector field} of $f$ by
\[
  V = - \frac{\nabla (\re f)}{| \nabla (\re f) |^2}
    = \frac{\xi_{\im f}}{|\xi_{\im f}|^2} .
\]
The flow of $V$ is called the gradient-Hamiltonian flow.
From the definition, we have 
\[
  V (\re f) = - \frac 1{| \nabla (\re f) |^2}
  \langle \nabla (\re f), \nabla (\re f) \rangle
  = -1
\] 
and
\[
  V (\im f) = \frac 1{|\xi_{\im f}|^2} \{\im f, \im f \}
  = 0.
\]
Therefore the gradient-Hamiltonian flow sends $X_1$ to another 
member $X_{1-t}$ of the family.
Note that $V$ does not preserve the symplectic structure 
$\widetilde\omega$,
because $V$ is normalized and hence not a Hamiltonian vector field.
However we can check that the restrictions 
$\omega_t := \widetilde{\omega} |_{X_t}$ to $X_t$ are preserved.
In other words, the gradient-Hamiltonian flow gives a map
\[
  \phi_t = \exp (tV) : 
  (X_1, \omega_1) \longrightarrow 
  (X_{1-t}, \omega_{1-t}).
\]
between symplectic varieties.

\begin{rem}[Ruan \cite{R2}]
If we write $f$ locally as $u/v$ for some holomorphic functions 
$u, v$, then $V$  can be written as
\begin{equation}
  V = \frac{-2 \re\, 
    \bigl( \bar{v}(\nabla v - t \nabla u) \bigr)}
           {|du - t dv |^2}.
  \label{vector}
\end{equation}
In particular, $V$ is smooth on the smooth part of $X_t$.
Note that (\ref{vector}) make sense even on the locus 
where $f$ is not defined, if it is smooth.
\end{rem}

\begin{prop}\label{prop:inv}
 Assume that $(\mathfrak{X}, \widetilde{\omega})$ has a
 Hamiltonian action of a compact Lie group $G$, which preserves 
 each $X_t$.
 Let $\mu : \mathfrak{X} \to \mathfrak{g}^*$ be the moment map.
 Then, for $h \in C^{\infty}(\mathfrak{g}^*)$, $\mu^* h$ is invariant 
 under the gradient-Hamiltonian flow of $f$.
\end{prop}

\begin{proof}
$G$-invariance of $f$ and Lemma \ref{Lem:Noether} implies that 
\[
  \xi_{\im f} (\mu^* h) = 0,
\]
which proves the proposition.
\end{proof}

We apply this to each stage $f_m : \mathfrak{X}_m 
\to \mathbb{C}$ of the toric degeneration.
We take a K\"ahler form on $\mathfrak{X}_m$, which is invariant
under the actions of 
$U(m-1)$ and $T_{n-1} \times \dots \times T_m$, 
and the restriction to each 
$X_{m,t}$ coincides with $\widetilde{\omega}_{\lambda}|_{X_{m,t}}$.
Then we obtain
the gradient-Hamiltonian vector field $V_m$ of $f_m$ and its flow 
\[
  \phi_{m,t} = \exp (tV_m) : 
  (X_{m,1}, \widetilde{\omega}_{\lambda}|_{X_{m,1}}) \longrightarrow 
  (X_{m,1-t}, \widetilde{\omega}_{\lambda}|_{X_{m,1-t}}).
\]
Using Proposition \ref{prop:inv} we have

\begin{cor}\label{Cor:flow}
  The gradient-Hamiltonian vector field $V_m$ preserves the 
  values of $\tilde{\nu}_i^{(n-1)}, \dots, \tilde{\nu}_j^{(m)},
  \tilde{\lambda}_k^{(m-1)}, \dots, \tilde{\lambda}_1^{(1)}$,
  and hence $\phi_{m,t}$ gives a deformation of $X_{m,t}$ 
  preserving the structure of completely integrable systems.
  In particular, the image $\Phi_{m,t} (X_{m,t})$ is 
  the Gelfand-Cetlin polytope $\Delta_{\lambda}$ for
  $t \ge 0$\rm{:}
  \[
    \xymatrix{
      X_{m,1}  \ar[dr]_{\Phi_{m,1}} \ar[rr]^{\phi_{m,1-t}}
      & & X_{m,t} \ar[dl]^{\Phi_{m,t}} \\ 
      & \Delta_{\lambda} &
    }
  \]
\end{cor}

Combining this with Theorem \ref{Thm:CIS}, we obtain a toric degeneration
of the Gelfand-Cetlin system.

\begin{rem}
 By changing the phase of the rational function $f_m$, we obtain a flow
 sending $X_{m, t}$ to $X_{m,0}$ for each $t$ not necessarily real.
 Note that, if $|t_0| = 1$, then 
 $(X_{n, t_0}, \widetilde{\omega}_{\lambda}|_{X_{n,t_0}})$
 is isomorphic to $(F, \omega_{\lambda})$, and 
 the restriction $\Phi_{n, t_0}$ coincides with the Gelfand-Cetlin
 system.
 Applying the same argument, we have 
 $\Phi_{m,t} (X_{m,t}) = \Delta_{\lambda}$ for each
 $t \in \mathbb{C}$.
\end{rem}

\begin{rem}
 The toric degeneration of the Gelfand-Cetlin system gives 
 an isomorphism between geometric quantizations for 
 the flag manifold and Gelfand-Cetlin toric variety.
 To see this, we recall the method of geometric quantization 
 via Lagrangian torus fibrations.
 Let $(M, \omega)$ be a symplectic manifold, $\mathcal{L}$ 
 a complex line bundle on $M$ with a unitary connection whose first
 Chern form coincides with $\omega$ (such a line bundle is called a
 {\it prequantum bundle}).
 We further assume that $M$ admits a Lagrangian torus fibration
 $\Phi : M \to B$.
 A fiber $L(u) := \Phi^{-1}(u)$ of $\Phi$ is said to be 
 {\it Bohr-Sommerfeld} if the restriction $\mathcal{L}|_{L(u)}$ has 
 trivial holonomies.
 The {\it real quantization} is defined to be the space of covariantly 
 constant sections of $\mathcal{L}$ restricted to Bohr-Sommerfeld fibers.
 
  Assume that $\lambda_{n_i}-\lambda_{n_{i+1}} \in \mathbb{Z}$
  for all $i$, and consider a family of line bundles
  $\mathfrak{L}_{\lambda} \to \mathfrak{X}_m$ given by
  \[
       \mathcal{O}_{\mathbb{P}_1}
     (\lambda_1 - \lambda_2) 
     \boxtimes \dots \boxtimes 
     \mathcal{O}_{\mathbb{P}_{n-1}}
     (\lambda_{n-1} - \lambda_n) ,
  \]
  on $\prod_i \mathbb{P}_{n_i}$.
  Then the restriction $\mathcal{L}_{m,t} = \mathfrak{L}_{\lambda}|_{X_{m,t}}$
  gives a prequantum line bundle on each $X_{m,t}$.
  Using the completely integrable systems $\Phi_{m,t}$, 
  we obtain a real quantization for each $X_{m,t}$.
  It is proved in \cite{GS} that Bohr-Sommerfeld fibers for the 
  Gelfand-Cetlin system exist exactly on integral points
  of the Gelfand-Cetlin polytope.
  It is easy to check that the vector field $V_m$ 
  lifts to $\mathfrak{L}_{\lambda}$
  preserving the unitary connection on each
  $\mathcal{L}_{m,t}$.
  In particular, the gradient-Hamiltonian flow preserves the 
  Bohr-Sommerfeld condition,
  and hence it gives an isomorphism between real quantizations
  on the flag manifold and the Gelfand-Cetlin toric variety.
\end{rem}

\begin{expl}
We see the gradient-Hamiltonian flow for $F^{(3)}$ in 
$\mathbb{P}^2 \times \mathbb{P}^2$
instead of in the total space $\mathfrak{X}$ of the deformation.
Recall that the degenerating family is given by
\[
  X_t =
  \Bigl\{ 
    \bigl([Z_1 : Z_2 : Z_3], [Z_{12} : Z_{13} : Z_{23}] \bigr)
   \, \Bigm| Z_1 Z_{23} - Z_2 Z_{13} + t Z_3 Z_{12} =0 \,
   \Bigr\}.
\]
Hence the rational function in this case is
\[
  f = \frac{Z_2 Z_{13} - Z_1 Z_{23}}{Z_3 Z_{12}}.
\]
Theorem \ref{Thm:CIS} says that the restriction 
\[
  \Phi_t := 
  (\tilde{\lambda}_1^{(2)}, \tilde{\lambda}_2^{(2)}, 
   \tilde{\lambda}_1^{(1)})|_{X_t} 
  : (X_t , \widetilde{\omega}_{\lambda}|_{X_t})
  \longrightarrow 
  \mathbb{R}^3
\]
to $X_t$ is a completely integrable system for each $t$, and 
$\Phi_0$ coincides with the moment map of the torus action 
on the Gelfand-Cetlin toric variety.
The gradient-Hamiltonian vector field $V$ vanishes on 
$X_1 \cap X_0$, which is the inverse image of two faces 
in the back side of $\Delta_{\lambda}$ in Figure \ref{GCpoly}.
Hence $X_1 \cap X_0$ is fixed under the gradient-Hamiltonian flow.
We see the behavior of the $S^3$-fiber of the Gelfand-Cetlin system.
Recall that this $S^3$ is given by
$\lambda_1^{(2)} = \lambda_2^{(2)}= \lambda_1^{(1)} = \lambda_2$.
It is easy to check that the image in $\mathbb{P}^2 \times \mathbb{P}^2$
of the $S^3$-fiber is
\[
  \Bigl\{ \bigl( [z_1 : z_2: \lambda_1- \lambda_2],
  [\lambda_2 - \lambda_3: \bar{z}_2: - \bar{z}_1] \bigr)
  \Bigm| |z_1|^2 + |z_2|^2 
  = (\lambda_1- \lambda_2)(\lambda_2 - \lambda_3) \Bigr\}.
\]
From Corollary \ref{Cor:flow}, the image of $S^3$
under the flow is given by 
$\Phi_t^{-1}(\lambda_2, \lambda_2, \lambda_2)$,
or equivalently,
\[
  \phi_{1-t} (S^3) = \left\{ \mu^{(2)} (Z) = 
  \begin{pmatrix} \lambda_2 & \\
                     & \lambda_2 \end{pmatrix} \right\}
  \cap X_t.
\]
Let $(y_1, y_2, y_{23}, y_{13})$ be local coordinates on
$\mathbb{P}^2 \times \mathbb{P}^2$ given by
\[
  y_i = \frac{Z_i}{\sqrt{\sum_j |Z_j|^2}}, \quad
  y_{ij} = \frac{Z_{ij}}{\sqrt{\sum_{k,l} |Z_{kl}|^2}}.
\]
Then we have
\begin{align*}
  \mu^{(2)} 
  &= 
  \frac{\lambda_1 - \lambda_2}{\sum |Z_i|^2}
  \begin{pmatrix}  |Z_1|^2 & \overline{Z}_1 Z_2 \\
                   Z_1 \overline{Z}_2 & |Z_2|^2
  \end{pmatrix} \\
  &\qquad +
  \frac{\lambda_2 - \lambda_3}{\sum |Z_{ij}|^2}
  \begin{pmatrix}  
   |Z_{12}|^2 + |Z_{13}|^2 & \overline{Z}_{13} Z_{23} \\
    Z_{13} \overline{Z}_{23} & |Z_{12}|^2 + |Z_{23}|^2
  \end{pmatrix} 
  + \begin{pmatrix} \lambda_3 & \\
                     & \lambda_3
  \end{pmatrix} \\
  &= 
  \frac{\lambda_1 - \lambda_2}{\sum |Z_i|^2}
  \begin{pmatrix}  |Z_1|^2 & \overline{Z}_1 Z_2 \\
                   Z_1 \overline{Z}_2 & |Z_2|^2
  \end{pmatrix} 
  +
  \frac{\lambda_2 - \lambda_3}{\sum |Z_{ij}|^2}
  \begin{pmatrix}  
   -|Z_{23}|^2 & \overline{Z}_{13} Z_{23} \\
    Z_{13} \overline{Z}_{23} & - |Z_{13}|^2
  \end{pmatrix} 
  + \begin{pmatrix} \lambda_2 & \\
                     & \lambda_2
  \end{pmatrix}\\
  &= 
  (\lambda_1 - \lambda_2)
  \begin{pmatrix}  |y_1|^2 & \overline{y}_1 y_2 \\
                   y_1 \overline{y}_2 & |y_2|^2
  \end{pmatrix} 
  +
  (\lambda_2 - \lambda_3)
  \begin{pmatrix}  
   -|y_{23}|^2 & \overline{y}_{13} y_{23} \\
    y_{13} \overline{y}_{23} & - |y_{13}|^2
  \end{pmatrix} 
  + \begin{pmatrix} \lambda_2 & \\
                     & \lambda_2
  \end{pmatrix}.
\end{align*}
It is easy from this to see that $\phi_{1-t} (S^3)$ is given by
\[
  \phi_{1-t} (S^3) = \Bigl\{ \sqrt t y 
  = \bigl(\sqrt t y_1, \sqrt t y_2, 
    \sqrt t y_{23}, \sqrt t y_{13} \bigr)
  \Bigm|  y \in S^3  \Bigr\},
\]
and this means that the $S^3$-fiber shrinks to the singular point 
of $X_0$ under the flow.
\end{expl}




\section{A not-in-stages degeneration of the Gelfand-Cetlin system}
\label{sc:direct degeneration}

First we summarize what we have obtained and then we will mention
 what we will need
 for the application of the toric degeneration
 to the Floer theory of flag manifolds, which is the content of the latter half
 of the paper. 

So far we considered degeneration of flag manifolds in stages
 and we saw
 that the gradient-Hamiltonian flow connects the two integrable system
 structures: Gelfand-Cetlin and toric.
The degeneration is parametrized by $\bC^{n-1}$.
The point $(0, \dots, 0)$ corresponds to the Gelfand-Cetlin toric variety.
The point $(1, \dots, 1)$ corresponds to the flag manifold
 embedded in the product of projective spaces 
 by the Pl\"ucker embedding.
This has a natural action of $U(n)$ (which can be extended to
 an action on the ambient multiple projective space),
 and from this we construct the 
 Gelfand-Cetlin system.

The gradient-Hamiltonian flow maps the fiber over
 $(1, \dots, 1)$ to the fiber over $(1, \dots, 1, 0)$, then to the fiber over
 $(1, \dots, 1, 0, 0)$ and so on, along a piecewise linear path on the
 base.
From the point of view of finding structures of integrable systems,
 the problem is that when we move from  $(1, \dots, 1)$, 
 we cannot see the natural $U(n)$ action any more.
More precisely, since we use only $U(n-1)$ action to construct Gelfand-Cetlin
 systems, true problem emerges after the first degeneration.
For example, there is only natural $U(n-2)$-action on the fiber
 over $(1, \dots, t, 0)$, $0\leq t < 1$.
So we cannot apply the construction of the Gelfand-Cetlin integrable system
 on these fibers.
The only structures of integrable systems we have on the fibers
 over $(1, \dots, t, 0)$, $0 < t < 1$ are the ones
 which are the push forward of the Gelfand-Cetlin system on 
 the fiber over  $(1, \dots, 1)$ by the gradient-Hamiltonian flow.
 
On $(1, \dots, 1, 0)$, there is another natural integrable structure, induced
  by the combination of the $U(n-2)$ action and the new action of the
  torus $T^{n-1}$
  (Theorem \ref{Thm:CIS}).
The good point is that this natural structure of an integrable system
 on the fiber over $(1, \dots, 1, 0)$ coincides 
 with the structure 
 of the integrable system induced by the push forward of the
 Gelfand-Cetlin system by the gradient-Hamiltonian flow
 (Corollary \ref{Cor:flow}).
The same procedure can be applied along the path from $(1, \dots, 1, 0)$ to
 $(1, \dots, 1, 0, 0)$ and so on, 
 giving the degeneration of the Gelfand-Cetlin integrable system on the flag manifold
 to the 
 Lagrangian fibration of the Gelfand-Cetlin toric variety. 

What we want to do from now is to compute some Floer theoretical quantity
 of the flag manifolds.
Namely, we want to compute the {\it potential function} of the Lagrangian torus fibers
 of the Gelfand-Cetlin system.
Such a computation was done by Cho and Oh \cite{CO} for toric manifolds.
Since the  Gelfand-Cetlin system degenerates to the toric integrable system,
 we want to use
 their result in our case, too.
This works for $F^{(3)}$, which degenerates directly (not in stages) to the 
 toric integrable system,
 but in general, we cannot directly apply Cho-Oh's computation.
 
The first problem is we have to care about the singularity
 of the toric variety (\cite{CO} deals with smooth toric manifolds),
 but this point is not very problematic for the calculation of the potential functions,
 for which we need to consider only those disks with Maslov index two.
We study this point in the next section.

The other problem, which is also related to the first,
 is that since the gradient-Hamiltonian flow does not 
 preserve the complex structure, 
 the moduli space of holomorphic disks with Lagrangian boundary 
 condition may change along the flow
 (note that in the degeneration in stages, the variety near the Gelfand-Cetlin toric
 variety is singular, and not the flag manifold. To reach to the flag manifold by
 chasing back the gradient-Hamiltonian flow, we have to go the long way 
 along the piecewise linear path, which will change the complex structure   
 widely).
 
This problem will be resolved when we 
 have a fiber preserving flow on the family over
 the segment $(t, \dots, t)$, $0\leq t\leq 1$, 
 which preserves the symplectic forms and the Lagrangian fibration structures.
Note that in this case the fiber over $(0, \dots, 0)$ is the Gelfand-Cetlin toric variety, 
 and the others are flag manifolds.
Let $\epsilon$ be a positive small number.
The fiber over $(\epsilon, \dots, \epsilon)$ has structure of
 an integrable system , which
 is the push-forward of (so identical to) the Gelfand-Cetlin system.
This integrable system will sufficiently resemble the toric integrable 
 structure on the special fiber.
On the other hand, the complex structures of the fiber over
 $(\epsilon, \dots, \epsilon)$ and of the special fiber
 are quite close, too
 (at least away from the singular locus).
Note that the complex structure on the fiber over $(\epsilon, \dots, \epsilon)$
 is not the push-forward by the flow, but the natural one as a submanifold of
 the multiple projective space.

So we want to construct such a flow and 
 this is what we will do in the rest of this section
  (in fact, we will construct a flow over some piecewise linear path, not on
  $(t, \dots, t)$, since it will suffice for application.
  Moreover, we can easily modify the construction to actually have a flow on 
  $(t, \dots, t)$.).
One may think that the gradient-Hamiltonian flow along
 the segment $(t, \dots, t)$, $0\leq t\leq 1$ will give such a flow, 
 but this need not be true.
The gradient-Hamiltonian flow will produce a structure of 
 an integrable system  on
 the fiber over $(0, \dots, 0)$ by push-forward,
  but it need not coincide with the toric integrable structure in general
 (the resulting integrable system may depend on the path on the base).
So we will give another construction.
Recall that we constructed a family of projective varieties parametrized
 by ${\bf t} = (t_2, \dots, t_n)\in\bC^{n-1}$ in Section 4.
We denote the total space by $\widetilde{\mathfrak X}\to \bC^{n-1}$.
We denote the fiber over $(t_2, \dots, t_n)$ by
 $X_{(t_2, \dots, t_n)}$ as before.
The result is the following.

\begin{prop}
There is a toric degeneration $(\tilde{f}, \gamma, \widetilde{\Phi}, \phi)$
 (in the sense of  Definition
 \ref{def:toric_degeneration})
 of the Gelfand-Cetlin system
 $(X, \omega, \Phi)$ with the following properties.
\begin{enumerate}
\item $\tilde{f}: \widetilde{\mathfrak X}\to \bC^{n-1}$ 
 is the $(n-1)$-parameter family
 constructed in Section \ref{sc:deg_of_GC}.
\item $\gamma: [0, 1]\to \bC^{n-1}$ is a piecewise linear path
with $\gamma(0) = (1, \dots, 1)$ and $\gamma(1) = (0, \dots, 0)$.
There is a positive small number $\epsilon$ such that 
$\gamma$ restricted to the interval $[1-\epsilon, 1]$ is given by
 $\gamma(t) = (1-t, \dots, 1-t)$.
\item $\phi_t: X = X_1\to X_{1-t}$
 is a diffeomorphism for $t\neq 1$ and is surjective for $t = 1$.
 $\phi_t$ coincides with the gradient-Hamiltonian flow for $t\in [0, 1-\epsilon]$,
 but (possibly) \emph{not} for $t\in [1-\epsilon, 1]$.
\item Let $S$ be the singular locus of $X_0$ (which is (complex) 
 codimension three, see the comment after Theorem 4.1).
$\phi_t$ is a diffeomorphism for all $t\in [0, 1]$
 when we restrict it to $\phi^{-1}_1(X_0\setminus S)$.
In particular, the flow is a diffeomorphism when we restrict it to 
 $\phi^{-1}_1(X_0\setminus D)$, where $D$ is the union of the toric divisors.
\item $\phi_1$ restricted to $\phi^{-1}_1(X_0\setminus D)$ preserves the
symplectic and Lagrangian torus fibration structures. 
\end{enumerate}
\end{prop}
\proof
Let $\delta$ be a positive small number.
Consider the following piecewise linear path on $\bC^{n-1}$,
 which approximates the path $\Gamma$ over which we constructed the degeneration of
 the Gelfand-Cetlin system in stages.
We start from $(1, \dots, 1)$ and go straight to
 $(1, \dots, 1, \delta)$.
Then we turn and go to $(1, \dots, 1, \delta^2, \delta)$.
Proceeding similarly, we go to 
 $(\delta^{n-1}, \delta^{n-2}, \dots, \delta)$ through the piecewise linear path.
Then finally we go to $(\delta^{n-1}, \delta^{n-1}, \dots, \delta^{n-1})$.

Recall that the gradient-Hamiltonian flow can be defined when
 we have a one parameter family of hypersurfaces in a K\"ahler manifold.
In our case, over the segment
 between $(1, \dots, 1, \delta^i, \delta^{i-1}, \dots, \delta)$ and
 $(1, \dots, 1, \delta^{i+1}, \delta^{i}, \dots, \delta)$,
 we take the union of the fibers over a holomorphic disk containing this
 segment as the ambient K\"ahler manifold (the K\"ahler structure is induced
 by the restriction of the K\"ahler structure on the product of
 the disk and the multiple projective manifold).
The base parameter gives a one parameter family of hypersurfaces,
 so we have a gradient-Hamiltonian flow along the path
 between $(1, \dots, 1, \delta^i, \delta^{i-1}, \dots, \delta)$ and
 $(1, \dots, 1, \delta^{i+1}, \delta^{i}, \dots, \delta)$ for each $i$ and
 similarly along the path between
 $(\delta^{n-1}, \delta^{n-2}, \dots, \delta)$ and $(\delta^{n-1},
  \delta^{n-1}, \dots, \delta^{n-1})$.

Since the path from
 $(1, \dots, 1)$ to $(\delta^{n-1}, \delta^{n-2}, \dots, \delta)$ approximates
 the path $\Gamma$ (degeneration in stages), the integrable system 
 on the fiber over  $(\delta^{n-1}, \delta^{n-2}, \dots, \delta)$ 
 defined by the push forward of the Gelfand-Cetlin system
 by the gradient-Hamiltonian flow 
 approximates the integrable system on the fiber over $(\delta^n, 0, \dots, 0)$
 (which is also defined as the push forward of the Gelfand-Cetlin system),
 at least away from the singular locus
 (when one wants to be more precise, one can say as follows.
Consider the path from $(\delta^{n-1}, \delta^{n-2}, \dots, \delta)$ to
 $(\delta^{n-1}, 0, \dots, 0)$ and the gradient-Hamiltonian flow along this path.
 Then we have two structures of integrable systems on  the fiber over
 $(\delta^{n-1}, 0, \dots, 0)$, pushing forward the Gelfand-Cetlin system
 along two paths. 
 One can deform the one integrable structure to the other 
 on a complement of
 some compact subset (of small measure) of the singular locus
 by a diffeomorphism
 of small norm, bounded by $O(\delta)$.).

Since the path between $(\delta^{n-1}, \delta^{n-2}, \dots, \delta)$
 and $(\delta^{n-1}, \delta^{n-1}, \dots, \delta^{n-1})$ is very short,
 one sees that the induced integrable system on the fiber over 
 $(\delta^{n-1}, \dots, \delta^{n-1})$ also approximates the one on
 the fiber over $(\delta^{n-1}, 0, \dots, 0)$.
On the other hand, the structure of the integrable system 
 on the fiber over $(\delta^{n-1}, 0, \dots, 0)$ approximates that of the
 Gelfand-Cetlin toric variety.
Consequently, we constructed an integrable system, which is canonically identified
 with the Gelfand-Cetlin system, on the fiber over $(\delta^{n-1}, \dots, \delta^{n-1})$, 
 with the property that it approximates
 well the Gelfand-Cetlin toric integrable system on the special fiber
 in the sense remarked above.
 
Now we deform $\delta$ to 0.
Then for each point $p$ on the (half closed) segment
 between $(\delta^{n-1}, \dots, \delta^{n-1})$ and $(0, \dots, 0)$ ($(0, \dots, 0)$ is not
 contained in the segment), we have a diffeomorphism from 
 the fiber over $(1, \dots, 1)$ to the fiber over $p$, 
 by the composition of the gradient-Hamiltonian flows along the path.
Since the path deforms continuously as we deform $\delta$, the diffeomorphism 
 also deforms continuously (with respect to the $C^{\infty}$-norm)),
 so this defines a flow along the segment
 between $(\delta^{n-1}, \dots, \delta^{n-1})$ and $(0, \dots, 0)$.
By construction (and Corollary \ref{Cor:flow})
 it is clear that this flow extends to  $(0, \dots, 0)$, and
 the push forward of the Gelfand-Cetlin system converges to the
 toric integrable system.
This gives the desired flow.

Take $\gamma$ to be the piecewise linear 
 path 
\[(1, \dots, 1)\to(1, \dots, 1, \delta)\to\cdots
 \to (\delta^{n-1}, \dots, \delta)\to (\delta^{n-1}, \dots, \delta^{n-1})\to
 (0, \dots, 0),\]
  with a suitable parametrization,
 and define the flow $\phi$ as above, along $\gamma$.
The claims of the proposition are obvious consequences of the above
 push-forward construction and the convergence property of the integrable
 system by the flow.\qed


\section{Holomorphic disks in a flag manifold} \label{sc:holo_disk}
Let $(\tilde{f}: \widetilde{\mathfrak X}
 \to \bC^{n-1}, \gamma, \widetilde{\Phi}, \phi)$
 be a toric degeneration of a Gelfand-Cetlin system
 as in Section \ref{sc:direct degeneration}.
Let $\epsilon$ be a positive small number.
We consider the subfamily over the closed segment
 between $(\epsilon, \dots, \epsilon)$ and $(0, \dots, 0)$
 with a flow constructed in the last section.
We write the fiber over $(t, \dots, t)$ as $X_t$.
We also fix a point $u \in \Int \Delta_\lambda$
in the interior of the Gelfand-Cetlin polytope and
write the Lagrangian fiber $\Phi_t^{-1}(u) \subset X_t$ as $L_t$.
Recall from Corollary \ref{Cor:flow}
that the flow
$$
 \phi_{t'} : X_t \to X_{t - t'}
$$
induces a diffeomorphism
$$
 \phi_{t'}|_{L_t} : L_t \simto L_{t - t'}.
$$

\begin{lem}
\label{Maslov-1}
The flow
$$
 \phi_\epsilon : X_{\epsilon} \to X_0
$$
induces an isomorphism
$$
 (\phi_\epsilon)_* : \pi_2(X_{\epsilon}) \simto \pi_2(X_0)
$$
of the homotopy groups.
\end{lem}
\begin{proof}
Let $S$ be the singular locus of $X_0$ and
$$
 \phi_{\epsilon}: X_{\epsilon}\to X_0
$$
be the map induced by the flow.
Note that the restriction
$$
 \phi_{\epsilon}|_{X_{\epsilon}\setminus \phi_{\epsilon}^{-1}(S)}
  : X_{\epsilon}\setminus \phi_{\epsilon}^{-1}(S) \to X_0 \setminus S
$$
is a diffeomorphism.
Recall that $\pi_2(X_0)$ is generated by torus-invariant curves.
Let 
$$
 p: \widetilde X_0\to X_0
$$
be a small resolution.
Since the fan for $\Xtilde_0$ is obtained from that for $X_0$
without adding one-dimensional cones,
for any torus-invariant curve $l$ in $X_0$,
there is a unique torus-invariant curve $\widetilde l$ in $\widetilde X_0$
mapped isomorphically to the curve $l$.
We think of $l$ and $\ltilde$ as inclusion maps.
Since $\widetilde X_0$ is nonsingular
and the exceptional locus has real codimension greater than two,
one can continuously move $\widetilde l$ to a map 
$$
 \widetilde l': S^2\to \widetilde X_0,
$$
so that the image is disjoint
from the exceptional locus of $\widetilde X_0$.
Then $l' = p\circ \widetilde l'$ is homotopic to $l$
(seen as a map by inclusion)
and it can be lifted to $X_{\epsilon}$ by $\phi_{\epsilon}^{-1}$. 
This proves that the map 
$$
 (\phi_{\epsilon})_*: \pi_2(X_{\epsilon})\to \pi_2(X_0)
$$
is surjective.

On the other hand, recall that $\pi_2(X_{\epsilon})$ is generated
by classes of rational curves (e.g., by Schubert subvarieties).
Let $\beta_1, \beta_2\in \pi_2(X_{\epsilon})$ be any classes
represented by (union of) rational curves.
Let
$$
 \varphi_{1, i}: C_{1, i}, \to X_{t_i}, \qquad
 \varphi_{2, i}: C_{2, i}, \to X_{t_i}
$$
where $t_i\to 0$ as $i\to \infty$, be two sequences of 
holomorphic maps from possibly disjoint unions of prestable rational curves,
representing the classes $\beta_1$ and $\beta_2$ respetively. 
These sequences can also be seen as sequences of holomorphic curves
in the fixed ambient space: 
$$
 C_i\to \prod\Bbb P_{n_j}.
$$
By Gromov compactness theorem,
we can assume there are limits $\varphi_1$ and $\varphi_2$
of $\varphi_{1, i}$ and $\varphi_{2, i}$ respectively.
Again by using the small resolution, 
we can deform $\varphi_1$ and $\varphi_2$ so that
their images do not intersect the singular locus of $X_0$,
and we can lift them to $X_t$.
Now assume that $\varphi_1$ and $\varphi_2$ give the same class
in $\pi_2(X_0)$.
Then the homotopy classes of their lifts in $X_t$ are also the same.
On the other hand,
these classes are $\beta_1$ and $\beta_2$ respectively
by construction.
Hence $\varphi_{1, i}$ and $\varphi_{2, i}$ must be in the same homotopy class,
and the injectivity of the map $(\phi_{\epsilon})_*$ is proved.
\end{proof}

Lemma \ref{Maslov-1} and the long exact sequences of homotopy groups
for the pairs $(X_\epsilon, L_\epsilon)$ and $(X_0, L_0)$
immediately gives the following:

\begin{cor}\label{relative}
The flow
$$
 \phi_\epsilon : X_{\epsilon} \to X_0
$$
induces an isomorphism
$$
 (\phi_\epsilon)_* : \pi_2(X_{\epsilon}, L_{\epsilon}) \simto \pi_2(X_0, L_0)
$$
of the relative homotopy groups.
\end{cor}

The Maslov index of holomorphic disks into $(X_t, L_t)$ is a homomorphism
$$
 \mu: \pi_2(X_t, L_t) \to \bZ.
$$
Although $X_0$ is a singular variety,
we can define the Maslov index of disks into $(X_0, L_0)$
by using the isomorphism in Corollary \ref{relative}.
This is a reasonable definition
if any holomorphic disk into $(X_0, L_0)$ can be deformed
to avoid the singular locus of $X_0$,
so that they can be lifted to a (not necessarily holomorphic) map
into $(X_t, L_t)$.
Proposition \ref{prop:deform} below shows that this is indeed the case.
To state it,
we recall the notion of toric transversality of holomorphic curves
in a toric variety.

\begin{definition}[{Nishinou and Siebert \cite[Definition 4.1]{NS}}]
A holomorphic curve in a toric variety $X$
is said to be {\em torically transverse}
if it is disjoint from all toric strata of codimension greater than one.
A stable map $\varphi : C \to X$ is {\em torically transverse}
if $\varphi(C) \subset X$ is torically transverse and
$\varphi^{-1}(\Int X) \subset C$ is dense.
Here $\Int X$ is the complement of the toric divisors of $X$.
\end{definition}

Note that a torically transverse map
does not intersect the singular locus of the toric variety.

\begin{proposition} \label{prop:deform}
Any disk $\varphi: (D^2, S^1) \to (X_0, L_0)$ can be deformed
into a holomorphic disk with the same boundary condition
which is torically transverse.
\end{proposition}

We need the following to prove Proposition \ref{prop:deform}:

\begin{thm}[{Cho and Oh \cite[Theorem 5.3]{CO}}] \label{Cho-Oh}
Let $L$ be a Lagrangian torus fiber
in a smooth projective toric variety
$$
 X_\Sigma = (\bC^r \setminus Z(\Sigma)) / K.
$$
Here
$\fanr$ is the number of one-dimenional cones
of a fan $\Sigma$,
the subset $Z(\Sigma) \subset \bC^{\fanr}$ is defined
by the Stanley-Reisner ideal, and
$K$ is the kernel of the map $(\bCx)^{\fanr} \to (\bCx)^N$
defined by one-dimenional cones in $\Sigma$.
Then any holomorphic map
$$
 \varphi : (D^2, \partial D^2) \to (X_{\Sigma}, L)
$$
can be lifted to a holomorphic map
$$
 \varphitilde :D^2 \to \bC^r \setminus Z(\Sigma)
$$
so that the homogeneous coordinate functions
$
 (z_1(\varphitilde), \cdots, z_r(\varphitilde))
$
are given by the Blaschke products with constant factors;
$$
 z_j(\varphitilde)
  = c_j \cdot
     \prod_{k=1}^{\mu_j}
      \frac{z - \alpha_{j,k}}{1 - \overline{\alpha}_{j,k} z},
$$
where $c_j \in \bCx$,
$\alpha_{j, k}\in \Int D^2$ and
$\mu_j$ is a non-negative integer for $j = 1, \cdots, r$.
Moreover,
the Maslov index of $\varphi$ is given by
$$
 \nu(\varphi) = 2 \sum_{j=1}^r \mu_j.
$$
\end{thm}

\noindent
{\em Proof of Proposition \ref{prop:deform}.}
Let $\widetilde X_0$ be a small resolution of $X_0$
and $\psi$ be the proper transform of $\varphi$.
Since $\widetilde X_0$ is smooth, 
the map $\psi$ has an explicit description
$$
 z_{j}(\widetilde{\psi})
  = c_j \cdot 
     \prod_{k=1}^{\mu_j}
      \frac{z-\alpha_{j,k}}{1-\overline{\alpha}_{j,k}z}
$$
by Theorem \ref{Cho-Oh}.
Note that ${\psi}$ intersects a toric stratum of higher codimension
exactly when there are $j_1\neq j_2$ such that 
$\alpha_{j_1, k_1} = \alpha_{j_2, k_2}$ for some $k_1$ and $k_2$. 
From this remark, and
since $\widetilde{X_0}$ is a small resolution of $X_0$ so that
the exceptional locus has codimension larger than one,
we can make $\psi$ torically transverse by perturbing $\alpha_{j, k}$.
Since the resolution is small,
torically transverse disks in $\widetilde{X_0}$ project to 
torically transverse disks in $X_0$.
This proves the proposition.
\qed

\begin{corollary} \label{cr:singular}
Assume that a holomorphic map
$
 \varphi: (D^2, S^1)\to (X_0, L_0)
$
intersects the singular locus of $X_0$.
Then the Maslov index of $\varphi$ is larger than two. 
\end{corollary}
\begin{proof}
In the proof of Proposition \ref{prop:deform},
a disk $\varphi$ intersecting the singular locus
lifts to a disk in $\widetilde{X_0}$ 
whose description via Theorem \ref{Cho-Oh}
has at least two non-constant factors.
Hence when we deform it into torically transverse disk,
it intersects the toric boundary at least at two points.
This implies that $\varphi$ has Maslov index larger than two.
\end{proof}

\begin{lemma} \label{lm:W}
There is a small neighbourhood $W_0$
of the singular locus $S \subset X_0$
such that any holomorphic disk
$
 \varphi: (D^2, S^1)\to (X_0, L_0)
$
of Maslov index two is disjoint from the closure of $W_0$.
\end{lemma}

\begin{proof}
Theorem \ref{Cho-Oh} and Corollary \ref{cr:singular} implies that
any holomorphic disk with Maslov index two
is given by
$$
 z_j(\varphitilde) =
 \begin{cases}
   \displaystyle{
    c_i \cdot \frac{z - \alpha}{1 - \overline{\alpha} z}
   } & \text{if } j = i, \\
   c_j & \text{otherwise},
 \end{cases}
$$
for some $i \in \{ 1, \dots, r \}$,
$(c_j)_{j=1}^r \in (\bCx)^r$ and $\alpha \in \Int D^2$.
The image of
$$
 [(c_1, \dots, c_{i-1}, 0, c_{i+1}, \dots, c_r)] \in X_0
$$
by the moment map is determined by the condition
$$
 \varphi(\partial D^2) \subset L_0,
$$
so that the image $\varphi(D^2)$
intersects the toric boundary of $X_0$
only at an interior point of a toric divisor,
whose image by the moment map is independent of
$(c_j)_j$ and $\alpha$,
and the lemma follows.
\end{proof}

Let us introduce the following notation:

\begin{definition}
Let $M$ be a K\"{a}hler manifold, $L$ be its Lagrangian submanifold, and
$\alpha \in \pi_2(M, L)$ be a relative homotopy class.
Then $\scMbar_1(M, L; \alpha)$ will denote the moduli space of
stable maps of degree $\alpha$
from a bordered Riemann surface of genus zero
with one marked point and with Lagrangian boundary condition.
The open subspace of $\scMbar_1(M, L; \alpha)$
consisting of maps from a disk will be denoted by
$\scM_1(M, L; \alpha)$.
\end{definition}

\begin{theorem}[{Cho and Oh \cite[Theorem 6.1]{CO}}]
 \label{th:Cho-Oh_regularity}
Let
$
 X_\Sigma = (\bC^r \setminus Z(\Sigma)) / K
$
be a projective toric variety
and $L \subset X_\Sigma$ be a Lagrangian torus fiber.
Assume that a holomorphic disk
$$
 \varphi : (D^2, S^1) \to (X_\Sigma, L)
$$
is disjoint from the singular locus of $X_\Sigma$
and admits a lift
$$
 \varphitilde : (D^2, S^1) \to (\bC^r \setminus Z(\Sigma), \pi^{-1}(L))
$$
to the homogeneous coordinate space.
Then $\varphi$
is Fredholm regular.
\end{theorem}

Theorem \ref{th:Cho-Oh_regularity} shows that
$\scM_1(X_0 \setminus W_0, L_0; \alpha)$ is a smooth manifold
without any virtual structure.
We also have:

\begin{lemma}
If $\alpha\in\pi_2(X_0, L_0)$ is a class with Maslov index two,
then the evaluation map induces a diffeomorphism
$$
 \ev : \scM_1(X_0 \setminus W_0, L_0; \alpha) \simto L_0.
$$
\end{lemma}

\begin{proof}
This is clear from Theorem \ref{Cho-Oh} and Lemma \ref{lm:W}.
\end{proof}

Lemma \ref{lm:W} also shows the following:

\begin{lemma} \label{lm:toric_surj2}
If $\alpha \in \pi_2(X_0, L_0)$ is a class with Maslov index two,
then the inclusion
$$
 \scM_1(X_0 \setminus W_0, L_0; \alpha)
  \hookrightarrow \scM_1(X_0, L_0; \alpha)
$$
is surjective.
\end{lemma}

The fact that $X_0$ is a Fano variety implies the following:

\begin{lemma} \label{lm:toric_surj1}
If $\alpha\in\pi_2(X_0, L_0)$ be a class with Maslov index two,
then the natural inclusions
$$
\scM_1(X_0, L_0; \alpha)
  \hookrightarrow \scMbar_1(X_0, L_0; \alpha)
$$
is surjective.
\end{lemma}

\begin{proof}
Let
$$
 \varphi : D_1 \cup \dots \cup D_p \cup S_1 \cup \dots \cup S_q
  \to X_0
$$
be a stable map of genus zero and Maslov index two,
where $D_i$ and $S_i$ are disk and sphere components
of the domain curve.
Then the contribution of each $D_i$ to the Maslov index of $\varphi$
is greater than one by Theorem \ref{Cho-Oh},
and that of each $S_i$ is positive since $X_0$ is a Fano variety.
This implies $p = 1$ and $q = 0$
so that the lemma follows.
\end{proof}

The following is the main result in this section:

\begin{proposition} \label{prop:disk}
For any relative homotopy class
$
 \alpha \in \pi_2(X_0, L_0) 
$
of Maslov index two,
there is a positive real number $0 < t \le 1$ and a diffeomorphism
$$
 \psi : \scMbar_1(X_0, L_0; \alpha) \to \scMbar_1(X_t, L_t; \alpha)
$$
such that the diagram
$$
\begin{CD}
 H_*(\scMbar_1(X_0, L_0; \alpha)) @>{\ev_*}>> H_*(L_0) \\
  @V{\psi_*}VV @VV{(\phi_t^{-1})_*}V \\
 H_*(\scMbar_1(X_t, L_t; \alpha)) @>{\ev_*}>> H_*(L_t) \\
\end{CD}
$$
is commutative.
\end{proposition}

The existence of the map $\psi$ comes from
the Fredholm regularity:

\begin{prop} \label{Maslov1}
For a class $\alpha\in\pi_2(X_0, L_0)$ with Maslov index two
and a sufficiently small positive number $t$,
there is a map
$$
 \psi : \scMbar_1(X_0, L_0; \alpha)
  \to \scMbar_1(X_t, L_t; \alpha)
$$
which is a diffeomorphism
into a connected component of $\scMbar_1(X_t, L_t; \alpha)$.
\end{prop}
\proof
Let $\varphi_0$ be an element of $\scMbar_1(X_0, L_0; \alpha)$.
Lemmas \ref{lm:toric_surj2} and \ref{lm:toric_surj1}
implies that $\varphi_0$ is a holomorphic map
$$
 \varphi_0 : (D^2, \partial D^2) \to (X_0 \setminus W_0, L_0),
$$
which is Fredholm regular
by Theorem \ref{th:Cho-Oh_regularity}.
Then for sufficiently small $t$,
the differential equation for holomorphic maps
$
 \varphi_t : (D^2, \partial D^2) \to (X_t, L_t)
$
near $\phi_t^{-1} \circ \varphi_0$ 
is a small perturbation of the equation
on $X_0$ which has the solution $\varphi_0$,
so this equation also has a solution and it is
Fredholm regular.
From this, it follows that
the moduli space $\scMbar_1(X_t, L_t; \alpha)$
contains a connected component diffeomorphic to 
$
 \scMbar_1(X_0, L_0; \alpha)
  = \scM_1(X_0 \setminus W_0, L_0; \alpha).
$
\qed

To show the surjectivity of $\psi$,
we use the following version of the Gromov compactness theorem:
\begin{thm}[{Ye \cite[Theorem 0.2]{Ye}}] \label{Ye}
Let $(M, \omega)$ be a compact symplectic manifold
and assume that we are given the following data:
\begin{itemize}
 \item
$
 \{ J_t \}_{t \in [0,1]}
$
is a smooth family of tame
almost complex structures on $M$,
 \item
$
 \{ N_t \}_{t \in [0,1]}
$
is a smooth family of compact totally real submanifolds,
 \item
$
 \{ t_i \}
$
is a strictly decreasing sequence in $[0, 1]$ converging to $0$, and
 \item
$
 \{ \varphi_{i} \}_{i \in \bN}
$
is a sequence of pseudo holomorphic disks in $(M, J_{t_i})$
with boundary on $N_{t_i}$.
\end{itemize}
Assume further that
the area of $\varphi_i$ is uniformly bounded by a positive constant.
Then there is a subsequence of $\{ \varphi_{t_i} \}_{i \in \bN}$
which converges to a stable $J_0$-holomorphic map
from a bordered Riemannian surface of genus $0$ in $M$
with boundary on $N_0$.
\end{thm}
Now we can prove the following:
\begin{cor}
\label{Maslov3}
For sufficiently small $t$
and a class $\alpha \in \pi_2(X_t, L_t)$ of Maslov index two,
one has an inclusion
$$
 \scM_1(X_t \setminus W_t, L_t; \alpha) \subset \Image \psi,
$$
where
$
 W_t = \phi_t^{-1}(W_0).
$
\end{cor}
\begin{proof}
Suppose that the statement is false.
Then there is a sequence $\{ t_i \}_{i \in \bN}$
converging to zero
and a sequence
$$
 \varphi_{i} : (D^2, S^1) \to (X_{t_i} \setminus W_{t_i}, L_{t_i})
$$
of holomorphic disks
not contained in $\Image \psi$.
By Theorem \ref{Ye},
we can assume that $\varphi_{i}$ converges to a stable map
$$
 \varphi: C \to X_0
$$
of Maslov index two
from a bordered Riemannian surface $C$ of genus 0.
Strictly speaking, we need to care about the singularity of $X_0$.
But one can argue as follows.
Note that $X_0$ is equivariantly embedded in the product of projective spaces
 with a natural torus action.
So the Lagrangian torus fiber $L_0$ extends to a Lagrangian torus
 $\widetilde L_0$ of the product of
 projective spaces.
It is easy to deform $\widetilde L_0$
 to totally real submanifolds $\widetilde L_t$ so that $L_t\subset \widetilde L_t$,
 since totally real condition is an open condition.
Now we can apply Theorem \ref{Ye}.

By Lemma \ref{lm:toric_surj1},
the stable map $\varphi$ is a holomorphic disk,
and Proposition \ref{Maslov1} implies that
$\varphi_{t_i}$ for sufficiently small $t_i$ are contained
in the family constructed there, a contradiction.
\end{proof}

\begin{lemma} \label{lm:flag_minimal_maslov}
For sufficiently small $t$,
the Maslov index of any holomorphic disk
$$
 \phi : (D^2, \partial D^2) \to (X_t, L_t)
$$
is greater than or equal to two.
\end{lemma}
\begin{proof}
Assume that the statement is false.
Then there is a sequence $\{ t_i \}_{i \in \bN}$
converging to zero
and a sequence
$$
 \varphi_{i} : (D^2, S^1) \to (X_{t_i}, L_{t_i})
$$
of holomorphic disks with Maslov index less than two.
Then as in the proof of Corollary \ref{Maslov3},
we will have a subsequence of $\varphi_{i}$
converging to a stable map
$$
 \varphi: C \to X_0
$$
of Maslov index less than two,
which contradicts the fact
that $X_0$ has no such stable maps.
\end{proof}

\begin{lemma} \label{lm:flag_surj1}
For sufficiently small $t$,
the natural inclusion
$$
 \scM_1(X_t, L_t; \alpha)
  \to \scMbar_1(X_t, L_t; \alpha)
$$
is surjective.
\end{lemma}
\begin{proof}
This follows from Lemma \ref{lm:flag_minimal_maslov}
and the fact that $X_t$ is Fano
in just the same way as in the proof of Lemma \ref{lm:toric_surj1}.
\end{proof}

\begin{lem}
\label{Maslov4}
For sufficiently small $t$
and a class $\alpha \in \pi_2(X_t, L_t)$ of Maslov index two,
the natural inclusion
$$
 \scM_1(X_t \setminus W_t, L_t; \alpha)
  \to \scMbar_1(X_t, L_t; \alpha)
$$
is surjective.
\end{lem}
\begin{proof}
Assume that the statement is false.
Then there is a sequence $\{ t_i \}_{i \in \bN}$
converging to zero
and a sequence
$$
 \varphi_{i} : (D^2, S^1) \to (X_{t_i}, L_{t_i})
$$
of holomorphic disks intersecting $W_{t_i}$.
As in the proof of Corollary \ref{Maslov3},
we can show that $\varphi_{t_i}$ converges to a holomorphic map
$$
 \varphi: D^2 \to X_0
$$
from a disk with Maslov index two
by taking a suitable subsequence if necessary.
Then the image of $\varphi$ must intersect the closure of $W_0$,
which contradicts our choice of $W_0$ in Lemma \ref{lm:W}.
\end{proof}

The commutativity of the diagram in Proposition \ref{prop:disk}
follows from the standard cobordism argument
on variations of moduli spaces under perturbations.

\section{Potential functions for Gelfand-Cetlin systems}
 \label{sc:potential}

In this section,
we recall the definition of the potential function
and compute it for Lagrangian torus fibers
of the Gelfand-Cetlin system.
Since our treatment here follows
Fukaya, Oh, Ohta and Ono \cite{FOOO_toric_I} closely,
we only give a sketch of the proof and refer the reader
to \cite{FOOO_toric_I} for further details.

Let
\begin{align*}
 \Lambda_0
  &= \left\{ \left.
      \sum_{i=1}^\infty a_i T^{\lambda_i} \, \right| \,
       a_i \in \bC, \ 
       \lambda_i \ge 0, \ 
       \lim_{i \to \infty} \lambda_i = \infty
     \right\} 
\end{align*}
be the Novikov ring
and
$$
\begin{array}{cccc}
 \frakv : & \Lambda_0 & \to & \bR \\
 & \rotatebox{90}{$\in$} & & \rotatebox{90}{$\in$} \\
 & \sum_{i=1}^\infty a_i T^{\lambda_i} & \mapsto
 & \min_{i} \{ \lambda_i \}_{i=1}^\infty
\end{array}
$$
be its valuation.
The maximal ideal and the quotient field of the local ring $\Lambda_0$
will be denoted by $\Lambda_+$ and $\Lambda$ respectively.

For a Lagrangian submanifold $L$ in a symplectic manifold $M$,
Lagrangian intersection Floer theory
equips the $\Lambda_0$-valued cochain complex of $L$
with the structure of an $A_\infty$-algebra
\cite{Fukaya_MHACFH, FOOO2006}.
By taking the canonical model,
one obtains an $A_\infty$-structure
$\{ \frakm_k \}_{k=0}^\infty$
on $H^*(L; \Lambda_0)$.
An element $b \in H^1(L; \Lambda_+)$
is called a {\em weak bounding cochain}
if it satisfies the {\em Maurer-Cartan equation}
\begin{equation} \label{eq:Maurer-Cartan}
 \sum_{k=0}^\infty \frakm_k(b, \dots, b)
  \equiv 0 \mod \PD([L]).
\end{equation}
The set of weak bounding cochains
will be denoted by $\scMwhat(L)$.
For any $b \in \scMwhat(L)$,
one can twist the Floer differential as
$$
 \frakm_1^b(x)
  = \sum_{k, l}
     \frakm_{k+l+1}
      (b^{\otimes k} \otimes x \otimes b^{\otimes l}).
$$
Maurer-Cartan equation implies
$\frakm_1^b \circ \frakm_1^b = 0$
and the resulting cohomology group
$$
 HF((L; b), (L; b))
  = \frac{\Ker(\frakm_1^b : H^*(L; \Lambda_0) \to H^*(L; \Lambda_0))}
         {\Image(\frakm_1^b : H^*(L; \Lambda_0) \to H^*(L; \Lambda_0))}
$$
will be called the {\em deformed Floer cohomology}.
The {\em potential function}
\begin{equation*}
 \po : \scMwhat(L) \to \Lambda_+
\end{equation*}
is defined by
\begin{equation} \label{eq:potential}
 \sum_{k=0}^\infty \frakm_k(b, \dots, b)
  = \po(b) \cdot \PD([L]).
\end{equation}

Now fix $\lambda$ as in \eqref{lambda} and
let $\Phi_\lambda : F = F(n_1, \dots, n_r, n) \to \Delta_\lambda$ be
the Gelfand-Cetlin system.
Let $v_i \in \bR^N$ be the primitive inward normal vector
of the $i$-th face of $\Delta_\lambda$
and choose $\tau_i \in \bR$
so that
$$
 \ell_i(u) = \langle v_i, u \rangle - \tau_i
$$
defines the $i$-th face of the Gelfand-Cetlin polytope $\Delta_\lambda$.
Here $\langle \bullet, \bullet \rangle$ is the standard inner product on $\bR^N$
as in Definition \ref{def:reflexive}.
The Lagrangian fiber
$
 \Phi^{-1}_{\lambda}(u)
$
over an interior point $u \in \Int \Delta_\lambda$
of the Gelfand-Cetlin polytope
will be denoted by $L(u)$.
We will identify $H^1(L(u); \Lambda_+)$
with $(\Lambda_+)^N$
using the angle coordinate
dual to the standard coordinate on the range
$\bR^N$ of the Gelfand-Cetlin system.
The following theorem
is a Gelfand-Cetlin analogue of
\cite[Proposition 3.2 and Theorem 3.4]{FOOO_toric_I}:

\begin{thm} \label{th:potential}
For any $u \in \Int \Delta_\lambda$,
one has an inclusion
$$
 H^1(L(u); \Lambda_+)
  \subset \scMwhat(L(u))
$$
and the potential function is given by
\begin{equation} \label{eq:toric_potential}
 \po^u(x)
  = \sum_{i=1}^\nofaces
     e^{\langle v_i, x \rangle} T^{\ell_i(u)}.
\end{equation}
\end{thm}

{\it Sketch of proof}.
Since the proof is completely parallel to
that of \cite[Proposition 3.2 and Theorem 3.4]{FOOO_toric_I},
we will be brief here
and refer the reader to \cite{FOOO_toric_I}
for further details.
Recall that $A_\infty$-structure
on the cochain complex of $L(u)$ is defined by
\begin{equation*}
 \frakm_k(a_1, \dots, a_k)
  = \sum_{\beta \in \pi_2(F, L(u))} \frakm_{k, \beta}(a_1, \dots, a_k),
\end{equation*}
\begin{equation*}
 \frakm_{k, \beta}(a_1, \dots, a_k)
  = (\ev_0)_!
     ([\scM_{k+1}(L(u), \beta)]^{\virt}
        \cap (\ev_1^* a_1 \cup \dots \cup \ev_k^* a_k))
     \cdot T^{\beta \cap \omega}
\end{equation*}
where $\scM_{k+1}(L(u), \beta)$ is the moduli space
of stable maps
with Lagrangian boundary condition
from a bordered Riemann surface of genus zero to $F$
with $k+1$ marked points on the boundary,
$[\scM_{k+1}(L(u), \beta)]^{\virt}$ is
its virtual fundamental chain,
and
\begin{equation*}
 \ev_i : \scM_{k+1}(L(u), \beta) \to L(u),
  \qquad i = 0, \dots, k
\end{equation*}
is the evaluation 
at the $i$-th marked point.
Since
$$
 \vdim \scM_{k+1}(L(u), \beta)
  = \dim L(u) + \mu(\beta) + k - 2,
$$
one has
$$
 \deg \frakm_{k, \beta}(b, \dots, b)
  = 2 - \mu(\beta)
$$
if $\deg b = 1$,
which is negative if $\mu(\beta) > 2$.
Hence only $\beta \in \pi_2(F, L(u))$ with $\mu(\beta) \le 2$
contribute
$
 \frakm_k(b, \dots, b) = \sum_\beta \frakm_{k, \beta}(b, \dots, b).
$
Since $F$ is simply-connected
and $\pi_2(L(u))$ is trivial,
the long exact sequence of relative homotopy groups gives
$$
 1 \to \pi_2 (F) \to \pi_2(F, L(u)) \to \pi_1(L(u)) \to 1,
$$
which splits since
$
 \pi_1 (L(u)) \cong \bZ^N
$
is a projective $\bZ$-module.
Let
$
 \beta_i \in \pi_2(F, L(u))
$
denote the lift of the primitive inward normal vector
$
 v_i \in \bZ^N \cong \pi_1(L(u))
$
of the $i$-th face $\partial_i \Delta_\lambda$
of the Gelfand-Cetlin polytope $\Delta_\lambda$,
which is represented by a disk
intersecting $\Phi_{\lambda}^{-1}(\partial_i \Delta_{\lambda})$
transversely at one point and
disjoint from $\Phi_{\lambda}^{-1}(\partial_j \Delta_{\lambda})$
for $i \ne j$.
The result of \S \ref{sc:holo_disk} shows that
one can choose a complex structure on $F$
so that
\begin{itemize}
 \item
any holomorphic disk of Maslov index two
with Lagrangian boundary condition
is Fredholm regular,
 \item
there is no holomorphic disk with Lagrangian boundary condition
of Maslov index less than two, and
 \item
$\scM_1(L(u), \beta_i) \cong L(u)$
for any face $i$
of the Gelfand-Cetlin polytope $\Delta_\lambda$,
and this is the only case where
$\scM_1(L(u), \beta)$ is non-empty for $\mu(\beta) = 2$.
\end{itemize}
Now we can follow the same reasoning as in
\cite[\S 10]{FOOO_toric_I}:
The fact that any $x \in H^1(L(u); \Lambda_+)$ is a weak bounding cochain
follows from the non-existence of holomorphic disks
with Lagrangian boundary condition
whose Maslov index is less than two.
Since $\scM_{k+1}(L(u), \beta_i)$ has an open dense subset
diffeomorphic to the product $C_k \times L(u)$
of $L(u)$ and the configuration space
\begin{equation}
 C_k = \{ (t_1, \dots, t_k) \mid 0 < t_1 < \dots < t_k < 1 \},
\end{equation}
one has
\begin{align*}
 \int_{L(u)} \frakm_{k, \beta_i}(x, \dots, x)
  &= \int_{C_k \times L(u)} (\ev_1^* x \cup \cdots \cup \ev_k^* x)
      T^{\beta_i \cap \omega} \\
  &= \Vol(C_k) \left( \int_{\beta_i} x \right)^k
      T^{\beta_i \cap \omega} \\
  &= \frac{1}{k!} {\langle v_i, x \rangle}^k T^{\ell_i(u)}
\end{align*}
and the potential function is given by
\begin{align*}
 \po^u(x)
  &= \int_{L(u)} \sum_{k=0}^\infty \frakm_k(x, \dots, x) \\
  &= \sum_{i=1}^\nofaces \sum_{k=0}^\infty
      \int_{L(u)} \frakm_{k, \beta_i}(x, \dots, x) \\
  &= \sum_{i=1}^\nofaces \sum_{k=0}^\infty
      \frac{1}{k!} {\langle v_i, x \rangle}^k T^{\ell_i(u)} \\
  &= \sum_{i=1}^\nofaces e^{\langle v_i, x \rangle} T^{\ell_i(u)}.
\end{align*}
This concludes the proof of Theorem \ref{th:potential}.
\qed

%
%

The following is an immediate corollary
of Theorem \ref{th:potential}:

\begin{cor}
The potential function
$$
 \po^u : H^1(L(u); \Lambda_+) \to \Lambda_+
$$
can be regarded as a Laurent polynomial
$$
 \po^u \in \bQ[Q_1^{\pm 1}, \dots, Q_{r+1}^{\pm 1}]
             [y_1^{\pm 1}, \dots, y_N^{\pm 1}]
$$
where
$$
 y_k = e^{x_k} T^{u_k}, \qquad k = 1, \dots, N
$$
are combinations of the variable
$x \in H^1(L(u); \Lambda_+)$
with the parameter $u \in \Delta_\lambda$
for the position of the fiber and
$$
 Q_j = T^{\lambda_{n_j}}, \qquad j = 1, \dots, r+1
$$
is the parameter
for the symplectic structure on $F$.
\end{cor}

\section{Examples} \label{sc:example}

In this section,
we study the critical points
of the potential function
for the full flag manifold $F(1, 2, 3)$
and the Grassmannian $\Gr(2, 4)$.
In the latter case,
we will see that the number of critical points
is strictly smaller
than the rank of the cohomology group.

Let us first discuss the case of $F(1, 2, 3)$.
Fix
$
 \lambda = (\lambda_1, \lambda_2, \lambda_3) \in \bR^3
$
satisfying
$$
 \lambda_1 > \lambda_2 > \lambda_3,
$$
so that the corresponding coadjoint (or adjoint) orbit
$\scO_{\lambda}$ is the full flag manifold of dimension three.
The Gelfand-Cetlin pattern in this case
is given by
\begin{equation*}
\begin{alignedat}{9}
 \lambda_1 &&&& \lambda_2 &&&& \lambda_3 \\
  & \uge && \dge && \uge && \dge \\
  && u_1 &&&& u_2 && \\
  &&& \uge && \dge &&& \\
  &&&& u_3 &&&&
\end{alignedat}
\end{equation*}
and the Gelfand-Cetlin polytope $\Delta_\lambda$
is defined by six inequalities
$$
 \Delta_\lambda
  = \{ u = (u_1, u_2, u_3) \in \bR^3 \mid
         \ell_i(u) = \langle v_i, u \rangle - \tau_i \ge 0,
          \quad i=1, \dots, 6
    \}
$$
where
\begin{align*}
 \ell_1(u) &= \langle (-1, 0, 0), u \rangle + \lambda_1, \\
 \ell_2(u) &= \langle (1, 0, 0), u \rangle - \lambda_2, \\
 \ell_3(u) &= \langle (0, -1, 0), u \rangle + \lambda_2, \\
 \ell_4(u) &= \langle (0, 1, 0), u \rangle - \lambda_3, \\
 \ell_5(u) &= \langle (1, 0, -1), u \rangle, \\
 \ell_6(u) &= \langle (0, -1, 1), u \rangle.
\end{align*}
The potential function is given by
\begin{align*}
 \po
  &= e^{-x_1} T^{- u_1 + \lambda_1}
     + e^{x_1} T^{u_1 - \lambda_2}
     + e^{- x_2} T^{- u_2 + \lambda_2} \\
  & \qquad
     + e^{x_2} T^{u_2 - \lambda_3}
     + e^{x_1 - x_3} T^{u_1 - u_3}
     + e^{- x_2 + x_3} T^{- u_2 + u_3} \\
  &= \frac{Q_1}{y_1} + \frac{y_1}{Q_2} + \frac{Q_2}{y_2}
     + \frac{y_2}{Q_3} + \frac{y_1}{y_3} + \frac{y_3}{y_2}.
\end{align*}
By equating the partial derivatives
\begin{align*}
 \frac{\partial \po}{\partial y_1}
  &= - \frac{Q_1}{y_1^2} + \frac{1}{Q_2} + \frac{1}{y_3}, \\
 \frac{\partial \po}{\partial y_2}
  &= - \frac{Q_2}{y_2^2} + \frac{1}{Q_3} - \frac{y_3}{y_2^2}, \\
 \frac{\partial \po}{\partial y_3}
  &= - \frac{y_1}{y_3^2} + \frac{1}{y_2}
\end{align*}
with zero,
one obtains
\begin{align*}
 Q_1 Q_2 y_3 &= y_1^2(y_3 + Q_2), \\
 Q_3( y_3 + Q_2) &= y_2^2, \\
 y_1 y_2 &= y_3^2,
\end{align*}
whose solutions are given by
\begin{align*}
 y_1 &= \frac{y_3^2}{y_2}, \\
 y_2 &= \pm \sqrt{Q_3( y_3 + Q_2)}, \\
 y_3 &= \sqrt[3]{Q_1 Q_2 Q_3}, \ 
        \omega \sqrt[3]{Q_1 Q_2 Q_3}, \ 
        \omega^2 \sqrt[3]{Q_1 Q_2 Q_3},
\end{align*}
where $\omega = \exp(2 \pi \sqrt{-1}/3)$
is a primitive cubic root of unity.
Since $\dim H^*(F(1, 2, 3), \Lambda)$ is six,
one has as many critical point as $\dim H^*(F(1, 2, 3), \Lambda)$
in this case.
One can show that
all these critical points are non-degenerate
by computing the Hessian.
The valuations of the critical points
are given by
\begin{align*}
 u_1
  &= \frakv (y_1)
   = \frakv (y_3^2 / y_2) \\
  &= - u_2 + 2 u_3, \\
 u_2
  &= \frakv (y_2)
   = \frac{1}{2} \frakv(Q_3 ( Q_2 + y_3)) \\
  &= \frac{1}{2}
      \left(
       \lambda_3 + \min\{ \lambda_2, u_3 \}
      \right), \\
 u_3
  &= \frakv (y_3)
  = \frac{1}{3} \frakv (Q_1 Q_2 Q_3) \\
  &= \frac{1}{3} ( \lambda_1 + \lambda_2 + \lambda_3),
\end{align*}
so that $u = (u_1, u_2, u_3)$ is unique
for any $\lambda = (\lambda_1, \lambda_2, \lambda_3)$
and always lies
in the interior of the Gelfand-Cetlin polytope.

Next we discuss the case of $\Gr(2, 4)$.
Fix
$
 \lambda = (\lambda_1, \lambda_2, \lambda_3, \lambda_4) \in \bR^3
$
satisfying
$$
 \lambda_1 = \lambda_2 > \lambda_3 = \lambda_4,
$$
so that $\scO_{\lambda}$ is the Grassmannian
of two-planes in a four-space.
The Gelfand-Cetlin pattern in this case is given by
\begin{equation*}
\begin{alignedat}{13}
  \lambda_1 &&&& \lambda_1 &&&& \lambda_3 &&&& \lambda_3 \\
  & \ueq && \deq && \uge && \dge && \ueq && \deq & \\
  && \lambda_1 &&&& u_1 &&&& \lambda_3 && \\
  &&& \uge && \dge && \uge && \dge &&& \\
  &&&& u_2 &&&& u_3 &&&& \\
  &&&&& \uge && \dge &&&&& \\
  &&&&&& u_4 &&&&&& \\
\end{alignedat}
\end{equation*}
so that the Gelfand-Cetlin polytope $\Delta_\lambda$
is defined by six inequalities
$$
 \Delta_\lambda
  = \{ u = (u_1, u_2, u_3, u_4) \in \bR^4 \mid
         \ell_i(u) = \langle v_i, u \rangle - \tau_i \ge 0,
          \quad i=1, \dots, 6
    \}
$$
where
\begin{align*}
 \ell_1(u) &= \langle (0, -1, 0, 0), u \rangle + \lambda_1, \\
 \ell_2(u) &= \langle (-1, 1, 0, 0), u \rangle, \\
 \ell_3(u) &= \langle (1, 0, -1, 0), u \rangle, \\
 \ell_4(u) &= \langle (0, 0, 1, 0), u \rangle - \lambda_3, \\
 \ell_5(u) &= \langle (0, 1, 0, -1), u \rangle, \\
 \ell_6(u) &= \langle (0, 0, -1, 1), u \rangle,
\end{align*}
and the potential function is given by
\begin{align*}
 \po
  &= e^{- x_2} T^{- u_2 + \lambda_1}
     + e^{- x_1 + x_2} T^{- u_1 + u_2}
     + e^{x_1 - x_3} T^{u_1 - u_3} \\
  & \qquad
     + e^{x_3} T^{u_3 - \lambda_3}
     + e^{x_2 - x_4} T^{u_2 - u_4}
     + e^{- x_3 + x_4} T^{- u_3 + u_4} \\
  &= \frac{Q_1}{y_2} + \frac{y_2}{y_1} + \frac{y_1}{y_3}
     + \frac{y_3}{Q_3} + \frac{y_2}{y_4} + \frac{y_4}{y_3}.
\end{align*}
By equating the partial derivatives
\begin{align*}
 \frac{\partial \po}{\partial y_1}
  &= - \frac{y_2}{y_1^2} + \frac{1}{y_3}, \\
 \frac{\partial \po}{\partial y_2}
  &= - \frac{Q_1}{y_2^2} + \frac{1}{y_1} + \frac{1}{y_4}, \\
 \frac{\partial \po}{\partial y_3}
  &= - \frac{y_1}{y_3^2} + \frac{1}{Q_3} - \frac{y_4}{y_3^2}, \\
 \frac{\partial \po}{\partial y_4}
  &= - \frac{y_2}{y_4^2} + \frac{1}{y_3}
\end{align*}
with zero,
one obtains
\begin{align*}
 y_1^2 &= y_2 y_3, \\
 Q_1 y_1 y_4 &= y_2^2 (y_1 + y_4), \\
 y_3^2 &= Q_3 ( y_1 + y_4), \\
 y_4^2 &= y_2 y_3,
\end{align*}
whose solutions are given by
\begin{align*}
 y_1 &= \pm \sqrt{Q_1 Q_3}, \\
 y_2 &= Q_1 Q_3 / y_3, \\
 y_3 &= \pm \sqrt{2 Q_3 y_1}, \\
 y_4 &= y_1.
\end{align*}
Since $\dim H^*(\Gr(2, 4), \Lambda)$ is six,
one has less critical point than $\dim H^*(\Gr(2, 4), \Lambda)$
in this case,
in contrast to the case of $F(1, 2, 3)$.
All these critical points and non-degenerate
and one can see that
\begin{align*}
 u_1
  &= \frac{1}{2} (\lambda_1 + \lambda_3), \\
 u_2
  &= \frac{1}{4} (3 \lambda_1 + \lambda_3), \\
 u_3
  &= \frac{1}{4} (u_1 + 3 \lambda_3), \\
 u_4
  &= \frac{1}{2} (\lambda_1 + \lambda_3),
\end{align*}
so that $u = (u_1, u_2, u_3, u_4)$ is unique
for any $\lambda = (\lambda_1, \lambda_2, \lambda_3, \lambda_4)$
and always lies
in the interior of the Gelfand-Cetlin polytope.


\section{Non-displaceable Lagrangian torus fibers}
 \label{sc:non-displaceable}

We give a proof of the following theorem
in this section:

\begin{thm} \label{th:non-displaceable}
Let 
$
 \lambda = (\lambda_1 \ge \lambda_2 \ge \cdots \ge \lambda_n)
$
be a non-increasing sequence of real numbers
and $\Phi_\lambda : F \to \Delta_{\lambda}$ be
the corresponding Gelfand-Cetlin system.
Then there exists $u \in \Int \Delta_\lambda$
such that the Lagrangian torus fiber
$
 L(u) = \Phi_{\lambda}^{-1}(u)
$
satisfies
$$
 \psi(L(u)) \cap L(u) \ne \emptyset
$$
for any Hamiltonian diffeomorphism
$
 \psi : \scO_\lambda \to \scO_\lambda.
$
If $\psi(L(u))$ is transversal to $L(u)$ in addition,
then
$$
 \#(\psi(L(u)) \cap L(u)) \ge 2^N.
$$
\end{thm}

This theorem is an analogue of \cite[Theorem 1.5]{FOOO_toric_I}
for flag manifolds,
and follows immediately
from Theorem \ref{th:floer_coh} below
and the Hamiltonian isotopy invariance
in Lagrangian intersection Floer theory
\cite[Theorem J]{FOOO2006}.

\begin{thm} \label{th:floer_coh}
For any $\lambda$,
there exists $u \in \Int \Delta_\lambda$
and $\frakx \in H^1(L(u); \Lambda_0)$
such that the deformed Floer cohomology is isomorphic
to the ordinary cohomology;
$$
 HF((L(u), \frakx), (L(u), \frakx) ; \Lambda_0)
  \cong H(L(u); \Lambda_0).
$$
\end{thm}

Note that $\frakx$ above is taken from $H^1(L(u); \Lambda_0)$
whereas the bounding cochain $b$
appearing in the definition of the deformed Floer cohomology
in section \ref{sc:potential}
is taken from $H^1(L(u); \Lambda_+)$.
To define the deformed Floer cohomology
twisted by $\frakx \in H^1(L(u); \Lambda_0)$,
one divide $\frakx$ into the constant part
and the positive part
\begin{align*}
 \frakx &= \frakx_0 + \frakx_+, \\
 \frakx_0 &\in H^1(L(u); \bC), \\
 \frakx_+ &\in H^1(L(u); \Lambda_+),
\end{align*}
take the flat non-unitary line bundle $\scL_{\rho}$
whose holonomy representation
is given by
$$
 \rho = \exp(\frakx_0) : H_1(L(u); \bZ) \to \bCx,
$$
consider the $A_\infty$-operation
$\{ \frakm_k^\rho \}_{k=0}^\infty$
twisted by the flat non-unitary line bundle $\scL_\rho$
as in Cho \cite{Cho},
and define the deformed Floer differential
$\frakm_1^\frakx$ by
$$
 \frakm_1^\frakx(x)
  = \sum_{k, l}
     \frakm_k^\rho(\frakx_+^{\otimes k} \otimes x
                    \otimes \frakx_+^{\otimes l}).
$$
Let $\{ \be_i \}_{i=1}^N$
be the basis of $H_1(L(u); \bZ)$
corresponding to the angle coordinate of $L(u)$.
Then for any $i = 1, \dots, N$, one has
\begin{align*}
 \frakm^{\frakx}_1(\be_i) \cap [L(u)]
  &= \sum_{k, l}
      \frakm^{\rho}_k
       (\frakx_+^{\otimes l} \otimes \be_i
         \otimes \frakx_+^{\otimes (k-l-1)})
       \cap [L(u)] \\
  &= \left. y_i\frac{\partial \po^u(y)}{\partial y_i} \right|_{y = \exp(\frakx)}.
\end{align*}
This shows that
$\frakm_1^{\frakx}= 0$ on $H^1(L(u); \Lambda_0)$
if $\fraky = \exp(\frakx)$ is a critical point of
$\po^u(y)$.
If this is the case,
the induction argument of \cite[Lemma 12.1]{FOOO_toric_I}
on the degree and the Maslov index $\mu(\beta)$
using the $A_\infty$-relation
\begin{align*}
 \frakm^{\rho, b}_{1, \beta}(\boldf_1 \cup \boldf_2)
  &= \sum_{\beta_1 + \beta_2 = \beta}
      \pm \frakm^{\rho, b}_{2, \beta_1}
           (\frakm^{\rho, b}_{1, \beta_2}(\boldf_1) \otimes \boldf_2) \\
  & \qquad + \sum_{\beta_1 + \beta_2 = \beta}
         \pm \frakm^{\rho, b}_{2, \beta_1}
              (\boldf_1 \otimes \frakm^{\rho, b}_{1, \beta_2}(\boldf_2)) \\
  & \qquad + \sum_{\beta_1 + \beta_2 = \beta}
         \pm \frakm^{\rho, b}_{1, \beta_1}
              (\frakm^{\rho, b}_{2, \beta_2}(\boldf_1 \otimes \boldf_2))
\end{align*}
shows that $\frakm^\frakx_1 = 0$ on $H^*(L(u); \Lambda_0)$,
so that the deformed Floer cohomology is isomorphic
to the ordinary cohomology;
$$
 HF^*((L(u), \frakx), (L(u), \frakx); \Lambda_0)
  \cong H^*(L(u); \Lambda_0).
$$
Note that
$$
 \frakv(\exp(\frakx)) = 0.
$$
Hence one can find a twisting cochain
$\frakx \in H^1(L(u); \Lambda_0)$
such that the deformed Floer cohomology
$
 HF^*((L(u), \frakx), (L(u), \frakx); \Lambda_0)
$
is isomorphic to the ordinary cohomology,
if there is a critical point
$$
 \fraky = (\fraky_1, \dots, \fraky_N) \in \Lambda^N
$$
of the Laurent polynomial
$$
 \frakP(y) = \sum_{i=1}^m y^{v_i} T^{- \tau_i}
$$
such that
$$
 \frakv(\fraky)
  = (\frakv(\fraky_1), \dots, \frakv(\fraky_N))
  \in \Delta_{\lambda}
  \subset \bR^N.
$$
The existence of such a critical point follows
from Proposition \ref{prop:valuation},
which we learned from Hiroshi Iritani.
See also \cite[Proposition 3.6]{FOOO_toric_I}.

\begin{prop} \label{prop:valuation}
For a convex polytope
$$
 \Delta
  = \{ u \in \bR^N
        \mid \ell_i(u) \ge 0, \ i = 1, \dots, m
    \}
$$
where
$$
 \ell_i(u) = \langle v_i, u \rangle - \tau_i,
$$
define a Laurent polynomial
$\frakP \in \Lambda[y_1^{\pm 1}, \dots, y_N^{\pm 1}]$ by
$$
 \frakP = \sum_{i=1}^m y^{v_i} T^{- \tau_i}.
$$
Then $\frakP$ has at least one critical point
whose valuation lies in the interior of $\Delta$.
\end{prop}


We divide the proof into three steps:

\begin{step} \label{st:real_positive_critpt}
Let
$
 P \in \bR[y_1^{\pm 1}, \dots, y_N^{\pm 1}]
$
be a Laurent polynomial over the field of real numbers
such that every non-zero coefficient is positive
and the origin is in the interior of the Newton polytope.
Then $P$ has a critical point
in $(\bR^{>0})^N$.
\end{step}

\begin{proof}
The set
$$
 \{ y \in (\bR^{>0})^N \mid P(y) \le c \}
$$
is compact for any $c \in \bR$,
so that $P$ has a global minimum in $(\bR^{>0})^N$.
\end{proof}

\begin{step} \label{st:Lambda_critpt}
The Laurent polynomial $\frakP$ has a critical point
$\fraky$ in $(\Lambda^{\times})^N$.
\end{step}

\begin{proof}
The set of critical points is defined as the common zero
of partial derivatives of the potential function,
which always exists in the compactification $\bP^N(\Lambda)$
of the torus $(\Lambda^{\times})^N$.
If all the critical points lie at infinity
and none of them lies on the torus,
then it remains so after substituting any real number into $T$.
However $\frakP$ has a critical point in $(\bR^{>0})^N$
after substituting any positive real number in $T$
by Step \ref{st:real_positive_critpt},
which shows that $\frakP$ also have a critical point
$\fraky$ on the torus $(\Lambda^{\times})^N$.
\end{proof}

\begin{step} 
The valuation of $\fraky$ lies
in the interior of $\Delta$.
\end{step}

\begin{proof}
Let $\Gamma$ be the convex hull of the set
$$
 \{ (v_i, z) \in \bR^{N} \times \bR \mid z \ge - \tau_i \}_{i}
   \cup \{ (0, z) \in \bR^{N} \times \bR \mid z \ge 0 \}
$$
and $\phi : \bR^N \to \bR$ be the piecewise-linear map
such that the union of faces of $\Gamma$
containing the origin is a part of the graph of $\phi$.
A subset of $\bR^N$ where $\phi$ is linear
forms a maximal-dimensional cone of a complete fan $\Sigma$ in $\bR^N$.
For each cone $\sigma$ in $\Sigma$,
define $u_\sigma \in \bR^N$ by
$$
 \phi|_{\sigma}(v) = - \langle u_\sigma, v \rangle.
$$
It follows from the construction of $\phi$ that
$$
 - \tau_i
  \ge \phi(v_i)
  \ge - \langle u_\sigma, v_i \rangle
$$
for any $i$, so that $u_\sigma \in \Delta$.
We will write $\Sigma^{(N)}$
for the set of $N$-dimensional cones
of the fan $\Sigma$.

Let $u$ be the valuation of $\fraky$
and put
$
 \tau_u = \min_{i} \{ \langle v_i, u \rangle - \tau_i \}.
$
The leading term $\frakP_u$ of $\frakP$ is defined as
$$
 \frakP_u
  = \sum_{i : \langle v_i, u \rangle - \tau_i = \tau_u}
     y^{v_i} T^{- \tau_i},
$$
which has the leading term $\fraky_0$ of $\fraky$
as its critical point.
Assume that $u$ is not in the interior of the convex hull of
$\{ u_\sigma \}_{\sigma \in \Sigma^{(N)}}$.
Then the Newton polytope of $\frakP_u$
will not contain the origin in its interior,
and one can choose a coordinate of the torus
so that $\frakP_u$ contains only non-negative powers of $y_1$.
This shows that the coefficient of any term in
$
 \partial \frakP_u / \partial y_1
$
is positive
if one substitutes a positive real number into $T$.
Recall from Step 2 that $\fraky$ gives positive real numbers
if one substitutes a positive real number into $T$.
It follows that $\fraky_0$ gives positive real numbers
after substituting a sufficiently small positive real number $\epsilon$
into $T$
and hence one has
$$
 \left. \frac{\partial \frakP_u}{\partial y_1}(\fraky_0)
  \right|_{T = \epsilon} > 0.
$$
This contradicts the fact that $\fraky_0$ is a critical point of $\frakP_u$
so that $u$ must be contained
in the interior of the convex hull of
$\{ u_\sigma \}_{\sigma \in \Sigma^{(N)}}$,
which in turn is contained in $\Delta$.
\end{proof}
\section{A relation with Toda lattice} \label{sc:Toda}

In this section,
we discuss the potential function
for the full flag manifold
after substituting $e^{-1}$
into the indeterminate element $T$ in the Novikov ring,
and its relation with quantum cohomology and the quantum Toda lattice.
Although the potential function is no longer invariant
under Hamiltonian isotopy
after this substitution
and hence unfit for application to symplectic topology,
it is the potential function after this substitution
which appears as the Landau-Ginzburg potential
of the mirror of Fano manifolds,
studied by string theorists
such as Hori and Vafa \cite{Hori-Vafa}.
The main result in this section is
Theorem \ref{th:GC_quantum_Toda},
which is an immediate consequence of Theorem \ref{th:potential}
and Givental's integral representation
in Theorem \ref{th:Givental_int_rep}.

%
%

Let us first recall the definition of quantum cohomology
and Givental's $J$-function.
For a projective manifold $X$
with its K\"{a}hler class $\omega$,
the quantum product $\circ$ on $H^*(X; \Lambda)$
is defined by
$$
 \langle A \circ B, C \rangle
  = \sum_{\beta \in H_2(X; \bZ)} T^{\beta \cap \omega}
     \int_{[\scMbar_{0, 3}(X, \beta)]^{\virt}}
      \ev_1^*(A) \cup \ev_2^*(B) \cup \ev_3^*(C)
$$
where
$
 \langle \bullet, \bullet \rangle
$
is the Poincar\'{e} pairing,
$
 [\scMbar_{0, 3}(X, \beta)]^{\virt}
$
is the virtual fundamental class
of the moduli space of stable maps
of genus zero and degree $\beta$
with three marked points into $X$,
and
$$
 \ev_i : \scMbar_{0, 3}(X, \beta) \to X, 
  \qquad i = 1, 2, 3,
$$
is the evaluation map at the $i$-th marked point.
The quantum product equips the cohomology group
$H^*(X; \Lambda)$
with the structure of a Frobenius algebra.
Now consider the substitution $T = e^{-1}$,
although it may not make sense
since the definition of quantum product involves an infinite sum.
When this sum converges,
the quantum cohomology ring $\circ$ can be regarded
as a family of Frobenius algebras,
parametrized by (an open subset of) $H^2(X; \bR)$
considered as the moduli space of symplectic structures.

Now choose a basis $\{ T_i \}_{i=1}^{\dimh}$
of $H^*(X; \bR)$
such that
$\{ T_i \}_{i=1}^r$ is a basis of $H^2(X; \bR)$.
Let $(t_i)_{i=1}^r$ be the coordinate of $H^2(X; \bR)$
dual to the basis $\{ T_i \}_{i=1}^r$,
so that the symplectic form $\omega$ of $X$ is represented as
$
 \omega = \sum_{i=1}^r t_i T_i.
$
Then quantum product
is an infinite series in
$$
 q = ( q_i )_{i=1}^r
   = ( \exp( - t_i) )_{i=1}^r,
$$
which,
in the case of the flag manifold,
is known to be convergent
for sufficiently small $q$.
One can also let $q$ take values
in the complexification $H^2(X; \bC)$
of $H^2(X; \bR)$.

The Givental's (small) $J$-function is defined by
$$
 J_j = \sum_{\beta \in H_2(F^{(n)}; \bZ)} q^{\beta}
  \int_{[\scMbar_{0, 1}(X); \beta)]^{\virt}}
   \frac{\ev^* (T_j \wedge \exp(\sum_{i=1}^{r} p_i t_i / \hbar))}
        {\hbar (\hbar - \psi)},
  \qquad j = 1, \dots, \dimh.
$$
It is known that
$$
 J_j = \langle s_j, 1 \rangle,
$$
where $(s_j)_{j=1}^\dimh$ is a basis of flat sections
of the {\em Givental connection},
which is a connection
on the trivial vector bundle
on $H^2(X; \bC)$
with fiber $H^*(X; \bC)$
defined by
$$
 \nabla_{\frac{\partial}{\partial t_i}}
  = \hbar \frac{\partial}{\partial t_i} - T_i \circ.
$$
The $\scD$-module on $H^2(X; \bC)$
generated by the $J$-function
is called the {\em quantum $\scD$-module},
whose characteristic variety is
the spectrum of the quantum cohomology ring.

For the full flag manifold $F^{(n)}$,
let $\scV_i \to F^{(n)}$ be the universal subbundle
of rank $i$ and
$$
 p_i
  = c_1(\scV_{i+1} / \scV_{i})
  \in H^2(F^{(n)} ; \bZ),
  \qquad i=0, \dots, n-1
$$
be the first Chern class
of the $i$-th quotient line bundle
$\scV_{i+1} / \scV_{i}$.
The set $\{ p_i \}_{i=0}^{n-1}$ generates $H^*(F^{(n)}; \bZ)$
and the complete set of relations is given by
$$
 (\lambda + p_0) \cdots (\lambda + p_{n-1})
  = \lambda^n.
$$
We introduce a redundant parameter $(t_i)_{i=0}^{n-1}$
for $H^2(F^{(n)}; \bC)$
and define the $J$-function by
$$
 J_j = \sum_{\beta \in H_2(F^{(n)}; \bZ)} q^{\beta}
  \int_{[\scMbar_{0, 1}(F^{(n)}); \beta)]^{\virt}}
   \frac{\ev^* (T_j \wedge \exp(\sum_{i=0}^{n-1} p_i t_i / \hbar))}
        {\hbar (\hbar - \psi)}
$$
where $j$ runs from $1$ to $\dim H^*(F^{(n)}; \bC) = (n-1)!$.

%
%

Now we recall
the quantum Toda lattice
and its relation with the $J$-function
of the full flag manifold
following Givental and Kim \cite{Givental-Kim}
(see also Kim \cite{Kim} and Joe and Kim \cite{Joe-Kim}).
The quantum Toda Hamiltonian is defined by
$$
 H = \frac{\hbar^2}{2}
      \sum_{i=0}^{n-1} \frac{\partial^2}{\partial t_i^2}
     - \sum_{i=1}^{n-1} e^{t_{i}- t_{i-1}}.
$$
It commutes with $n$ mutually commutative differential operators
\begin{equation} \label{eq:quantum_Toda}
 D_i \left(
      \hbar \frac{\partial}{\partial t_0}, \dots,
       \hbar \frac{\partial}{\partial t_{n-1}},
       q_1, \dots, q_{n-1}
     \right), \qquad
  i = 1, \dots, n
\end{equation}
where $q_i = \exp(t_i - t_{i-1})$,
$$
 \det(A + x I)
  = x^{n+1}
     + \sum_{i=1}^n D_i(p_0, \dots, p_{n-1}, q_1, \dots, q_{n-1}) x^{n-i},
$$
and
$$
 A = \begin{pmatrix}
      p_0 & q_1 &  0  & \cdots & 0 & 0 \\
       -1 & p_1 & q_2 & \cdots & 0 & 0 \\
       0  &  -1 & p_2 & \cdots & 0 & 0 \\
   \vdots & \vdots & \vdots & \ddots & \vdots & \vdots \\
       0  &  0  &  0  & \cdots & p_{n-2} & q_{n-1} \\
       0  &  0  &  0  & \cdots & -1 & p_{n-1}
     \end{pmatrix}.
$$

%
%

The following theorem gives an astonishing relation
between the quantum cohomology ring
of the full flag manifold
and the quantum Toda lattice:

\begin{thm}[{Givental and Kim \cite{Givental-Kim}, Kim \cite{Kim}}]
The $J$-function of the full flag manifold $F^{(n)}$
is an eigenfunction of the quantum Toda lattice:
$$
 D_i J_j = 0,
  \qquad i = 1, \dots, n,
  \quad  j = 1, \dots, n!.
$$
\end{thm}

It follows from this theorem that
the quantum cohomology ring
of the full flag manifold
is isomorphic to the coordinate ring
of the Lagrangian level set
of the classical Toda Hamiltonians.

%
%

Now we recall the stationary-phase integral representation
of the eigenfunction of the quantum Toda lattice
due to Givental \cite{Givental_SPI}.
Consider $n(n-1)$ variables
\begin{align*}
 \{ X_{ij}, Y_{ij} \mid i = 1, \dots, n-1, \ j = 1, \dots, n-i, \}
\end{align*}
and the $n(n-1)/2$-dimensional torus $Y_q$
cut out from
$
 \Spec \bC[X_{ij}^{\pm 1}, Y_{ij}^{\pm 1}]_{i,j}
$
by the equations
$$
 Y_{i, j} X_{i, j} = X_{i+1, j} Y_{i, j+1},
  \qquad i = 1, \dots, n - 2, \quad j = 1, \dots, n - i- 1,
$$
and
$$
 X_{i, n-i} Y_{i, n-i} = q_i,
  \qquad i = 1, \dots, n-1,
$$
where $q = (q_1, \dots, q_{n-1}) \in (\bCx)^{n-1}$.
These relations imply that
$X_{ij}$ and $Y_{ij}$ can be expressed
by $n(n-1)/2$ variables
$$
 \{ T_{ij} \mid i = 1, \dots, n - 1, \ j = 1, \dots, n - i \}
$$
as
$$
 X_{ij} = \exp(T_{ij} - T_{i, j+1})
$$
and
$$
 Y_{ij} = \exp(T_{i+1,j} - T_{ij}),
$$
where
$$
 q_i = \exp(T_{i+1, n - i} - T_{i, n-i+1}),
  \qquad i = 1, \dots, n-1.
$$
Define the phase function $f_q$
and the holomorphic volume form $\omega$
on $Y_q$ by
$$
 f_q = \sum_{i, j} (X_{ij} + Y_{ij})
$$
and
$$
 \omega = \bigwedge_{i,j} \  d T_{ij}.
$$
Fix $\hbar \in \bCx$
and a complete K\"{a}hler metric on $Y_q$.
A {\em Lefschetz thimble}
is the unstable manifold of $\Re (f_q / \hbar)$
starting from a critical point of $f_q$.
The Following theorem is due to Givental:
\begin{thm}[{Givental \cite{Givental_SPI}}]
 \label{th:Givental_int_rep}
The phase function has
$\dim H^*(F^{(n)}) = n!$
critical points,
and the stationary-phase integrals
$$ 
 I_a
  = \int_{\Gamma_a}
     e^{f_q/\hbar} \omega 
$$
for the corresponding Lefschetz thimbles
$\{ \Gamma_a \}_{a=1}^{n!}$
gives the component $J_a$ of the $J$-function
for a suitable choice
of a basis of $H^*(F^{(n)}; \bC)$.
\end{thm}

Now it is obvious that
the potential function $\po|_{T = e^{-1}}$
and the phase function $f_q$
are related by
$$
 T_{ij} = x^{(i+j-1)}_i + \lambda^{(i+j-1)}_i,
  \qquad i=1, \dots, n - 1,
  \quad j = 1, \dots, n - i,
$$
and
$$
 T_{i, n - i + 1} = \lambda_i
  \qquad i = 1, \dots, n.
$$
This results in the following
striking relation
between the Gelfand-Cetlin system
and the quantum Toda lattice:

\begin{thm} \label{th:GC_quantum_Toda}
The potential function for Lagrangian torus fibers
of the classical Gelfand-Cetlin system
on the full flag manifold $F^{(n)}$,
considered as a Laurent polynomial in $n(n-1)/2$ variables
with $n$ parameters after substituting $e^{-1}$ to $T$,
is the phase function for an integral representation
of the solution to the quantum Toda lattice.
\end{thm}


Now let us discuss the classical limit of the above story.
The classical Toda lattice
is a completely integrable system
whose Hamiltonians are the classical limits
$$
 D_i (p_0, \dots, p_{n-1}, q_1, \dots, q_{n-1}),
  \qquad i = 0, \dots, n-1
$$
of the differential operators
\eqref{eq:quantum_Toda}.
The level set of $\{ D_i \}_{i=0}^{n-1}$
is a Lagrangian subvariety of
$(\Spec \bC[p_1, \dots, p_{n-1}, q_1, \dots, q_{n-1}], \omega)$,
where $p_0$ is determined by
$$
 D_{1}(p_0, \dots, p_{n-1}) = p_0 + \dots + p_{n-1} = 0,
$$
and the symplectic form is given by
$$
 \omega = \sum_{i=1}^{n-1} p_i \wedge \frac{d q_i}{q_i}.
$$

The classical limit of
the stationary-phase integral
is controlled by the Jacobi ring
$$
 J(f_q) =
 \bQ[q_1, \dots, q_{n-1}][y_1^{\pm 1}, \dots, y_N^{\pm 1}]
      \left/ \left( \frac{\partial f_q}{y_i \partial y_i} \right)_{i=1}^N 
      \right.
$$
whose spectrum is
the set $\Cr(f_q)$ of critical points of $f_q$,
in that there is a birational map
$$
\begin{array}{ccc}
 \Cr(f_q) & \to & \Spec \bC[p_1, \dots, p_{n-1}, q_1, \dots, q_{n-1}] \\
 \rotatebox{90}{$\in$} & & \rotatebox{90}{$\in$} \\
 (y, q) & \mapsto &
  \displaystyle{
    \left( q \frac{\partial f_q}{\partial q}(y), \  q
    \right)
  },
\end{array}
$$
into the characteristic variety
of the $D$-module
generated by the stationary phase integrals.
On the other hand,
the characteristic variety of
the quantum $D$-module
is the spectrum of the quantum cohomology ring.
By putting them together,
we obtain the following:

\begin{cor}
The Jacobi ring of the potential function
for Lagrangian torus fibers
of the classical Gelfand-Cetlin system
on the full flag manifold $F^{(n)}$
is isomorphic to the ring of functions
on the level set of the classical Toda Hamiltonians,
and hence to the quantum cohomology ring of $F^{(n)}$.
\end{cor}

Note that the isomorphism
between the Jacobi ring of the potential function
for Lagrangian torus fibers
of the classical Gelfand-Cetlin system
and the quantum cohomology ring
can not hold
for general partial flag manifolds;
the simplest example is the Grassmannian $\Gr(2, 4)$
where the number of critical points
of the potential function
for general $q$ is four,
which is strictly smaller than the rank of the cohomology ring.


\bibliographystyle{plain}
\bibliography{gcbibs}


\noindent
Takeo Nishinou

Mathematical Institute,
Tohoku University,
Sendai,
980-8578,
Japan

{\em e-mail address}\ : \  nishinou@math.tohoku.ac.jp

\ \\

\noindent
Yuichi Nohara

Mathematical Institute,
Tohoku University,
Sendai,
980-8578,
Japan

{\em e-mail address}\ : \  nohara@math.tohoku.ac.jp

\ \\

\noindent
Kazushi Ueda

Department of Mathematics,
Graduate School of Science,
Osaka University,
Machikaneyama 1-1,
Toyonaka,
Osaka,
560-0043,
Japan.

{\em e-mail address}\ : \  kazushi@math.sci.osaka-u.ac.jp
\ \vspace{0mm} \\

Mathematical Institute,
University of Oxford,
24-29 St Giles'
Oxford
OX1 3LB

{\em e-mail address}\ : \  uedak@maths.ox.ac.uk

\end{document}